\documentclass[a4paper]{amsart}

\usepackage[totalwidth=17cm,totalheight=21cm]{geometry}
\usepackage{amssymb,amsmath,latexsym,pstricks}
\usepackage{diagrams}
\usepackage{amscd,epsfig}
\everymath{\displaystyle}

\newtheorem{teo}{Theorem}[section]
\newtheorem{lem}[teo]{Lemma}

\newtheorem{prop}[teo]{Proposition}
\newtheorem{defi}[teo]{Definition}

\newtheorem{remark}[teo]{Remark}
\newtheorem{remarks}[teo]{Remarks}

\newcommand{\mr}{\mathbb{R}}
\newcommand{\mc}{\mathbb{C}}
\newcommand{\mz}{\mathbb{Z}}
\newcommand{\mh}{\mathbb{H}}

\newcommand{\mpi}{\mathbb R P}

\newcommand{\Ee}{{\mathcal E}}

\newcommand{\Gg}{{\mathcal G}}
\newcommand{\Hh}{{\mathcal H}}

\newcommand{\Kk}{{\mathcal K}}
\newcommand{\Ll}{{\mathcal L}}
\newcommand{\Mm}{{\mathcal M}}

\newcommand{\Ss}{{\mathcal S}}
\newcommand{\Tt}{{\mathcal T}}

\newcommand{\R}{{\mathbb R}}
\newcommand{\Z}{{\mathbb Z}}
\newcommand{\cML}{\mathcal{ML}}
\newcommand{\cT}{{\mathcal T}}
\newcommand{\cTh}{\hat{\mathcal T}}

\newcommand{\cb}{\overline{c}}

\newcommand{\dr}{\partial}
\newcommand{\epsilonb}{\overline{\epsilon}}
\newcommand{\gammab}{\overline{\gamma}}
\newcommand{\lambdab}{\overline{\lambda}}
\newcommand{\Sigmab}{\overline{\Sigma}}
\newcommand{\St}{\tilde{S}}

\newcommand{\pf}{\begin{proof}}
\newcommand{\cvd}{\end{proof}}

\def\eps{\varepsilon}

\begin{document}

\title{Multi Black Holes and Earthquakes on Riemann surfaces with boundaries}
\author{Francesco Bonsante}
\address{Dipartimento di Matematica dell'Universit\`a degli Studi di Pavia,
via Ferrata 1, 27100 Pavia (ITALY)}
\email{francesco.bonsante@unipv.it}
\thanks{F. B. was partially supported by CNRS, ANR GEODYCOS.} 
\author{Kirill Krasnov}
\address{School of Mathematical Sciences, University of Nottingham, Nottingham,
  NG7 2RD, UK.}  
\email{kirill.krasnov@nottingham.ac.uk}
\thanks{K. K. was supported by an EPSRC Advanced Fellowship} 
\author{Jean-Marc Schlenker}
\address{Institut de Math\'ematiques de Toulouse (UMR CNRS 5219) \\
Universit\'e Paul Sabatier\\
31062 Toulouse cedex 9, France}
\email{jmschlenker@gmail.com}
\thanks{J.-M. S. was partially supported by the A.N.R. programs RepSurf,
2006-09, ANR-06-BLAN-0311, GeomEinstein, 2006-09, 06-BLAN-0154, and
 ETTT (NT09-512070).}

\date{July 2009 (v3)}

\begin{abstract}
We prove an ``Earthquake Theorem'' for hyperbolic metrics with geodesic
boundary on a compact surfaces $S$ with boundary:  given two hyperbolic metrics
with geodesic boundary on a surface with $k$ boundary components, there are
$2^k$ right earthquakes transforming the first in the second.  An alternative formulation 
arises by introducing the enhanced Teichm\"uller space of S: We prove that any two points 
of the latter are related by a unique right earthquake. The proof rests on the geometry 
of ``multi-black holes'', which are 3-dimensional anti-de Sitter manifolds, topologically 
the product of a surface with boundary by an interval.
\end{abstract}

\maketitle

\section{Introduction}

\subsubsection*{The Earthquake theorem}

Let $\Sigma$ be a closed surface, with a hyperbolic metric $g$, 
let $c$ be a simple closed geodesic on $(\Sigma, g)$, and let $l$
be a positive real number. The image of $g$ by the {\it right earthquake}
of length $l$ along $c$ is the hyperbolic metric obtained by cutting
$\Sigma$ along $c$ and gluing back after rotating the ``left'' side of
$c$ by $l$. This defines a map from the Teichm\"uller space $\cT_\Sigma$
of $\Sigma$ to itself. 

Suppose now that $\lambda$ is a measured geodesic lamination on $(\Sigma, h)$
which is {\it rational}, i.e., its support is a disjoint union of closed
curves $c_1, \cdots, c_n$. The transverse measure is then described by a 
set of positive numbers $l_1, \cdots, l_n$ associated to the $c_i$. 
The image of $g$ by the right earthquake along $\lambda$ is obtained
as above, by doing a ``fractional Dehn twist'' along each of the $c_i$,
with a length parameter given by the $l_i$. Again this defines a map
from $\cT_\Sigma$ to itself.

Thurston \cite{thurston-notes,thurston-earthquakes} 
discovered that this definition can be
extended by continuity to all measured geodesic laminations on $(\Sigma, g)$.
In other terms, it makes sense to talk about the right earthquake along
any measured geodesic lamination on $(\Sigma, g)$. This defines a map:
$$ E_r :\cML_\Sigma \times \cT_\Sigma \rightarrow \cT_\Sigma~, $$
where $\cML_\Sigma$ is the space of measured laminations on $\Sigma$. 
Thurston also discovered a striking feature of this Earthquake map.

\begin{teo}[Thurston \cite{thurston-earthquakes,kerckhoff}]\label{classical:teo}
For any $h,h'\in \cT_\Sigma$ there exists a unique $\lambda \in \cML_\Sigma$ 
such that $E_r(\lambda)(h)=h'$.
\end{teo}

\subsubsection*{Earthquakes on surfaces with boundary.}

Let now $\Sigma$ be a compact orientable surface of genus $g$ with $n$ boundary
components. We will assume $\Sigma$ to have
negative Euler characteristic
\[
   \chi(\Sigma)= 2-2g-n<0\,.
\]

Let $\cT_{g,n}$ be the Teichm\"uller
space of hyperbolic metrics on $\Sigma$ with geodesic boundary (such that each
geodesic boundary component is a closed curve), considered up
to isotopy. $\cT_{g,n}$ is a contractible manifold of dimension $6g-6+3n$. 

\begin{figure}
\begin{center}
\epsfig{figure=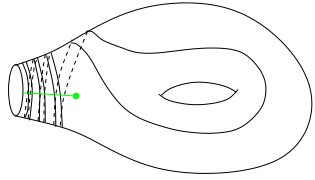,height=0.9in}
\end{center}
\caption{An example of a geodesic lamination on a surface with a geodesic
  boundary. The geodesics forming 
the lamination can spiral onto the boundary. The total weight of an arc ending
at the boundary (as shown in green here) is allowed 
to be infinite.}
\label{fig:lamin-open}
\end{figure}

We also consider the space $\cML_{g,n}$ of measured laminations on the 
{\it interior} of $\Sigma$, see e.g. \cite{FLP} (a precise definition
is given in section 3). Note that the transverse
weight on those laminations is required to be finite on any close transverse
segment in the interior of $\Sigma$, but the weight might be infinite on
segments with an endpoint on the boundary of $\Sigma$, see Figure \ref{fig:lamin-open}.
Given a measured lamination $\lambda\in \cML_{g,n}$ and a hyperbolic metric
$h\in \cT_{g,n}$, there is a unique way to realize $\lambda$ as a measured
{\it geodesic} lamination on $(\Sigma,h)$. 

The main result presented here is the following. 

\begin{teo} \label{tm:earthquake}
Given $h_1, h_2\in \cT_{g,n}$, there
are exactly $2^n$ measured laminations $\lambda_1, \cdots, \lambda_{2^n}$ 
on the interior of $\Sigma$ such that the right earthquake along the
$\lambda_i$ sends $h_1$ to $h_2$.
\end{teo}

This result extends to the hyperbolic metrics with some geodesic 
boundary components
and some cusps, however the number of possible measured laminations is lower
when one of the boundary components corresponds to a cusp for either $h_1$ or
$h_2$. The statement of Theorem \ref{tm:earthquake} looks simple, but it
might be less obvious than it first seems; even the case $g=0,n=3$ (for a hyperbolic
pair of pants), where everything can be described explicitly, displays some
interesting phenomena. This case is described in details at the end of section 
2 (see Proposition \ref{pr:pants-lamin} and the paragraph right before section 
3).

\subsubsection*{The enhanced Teichm\"uller space.}

The fact that the number of right earthquakes sending a given hyperbolic
metric to another one is $2^n$ rather than one can appear distressing at first
sight. There is a simple geometric formalism, however, under which this
disagreement disappears. It is based on a definition due to V. Fock 
\cite{fock97,fock-tchekhov1,fock-tchekhov2,bonahon-liu} which appeared
naturally in different contexts. The terminology is borrowed from Bonahon
and Liu \cite{bonahon-liu}.

\begin{defi} \label{df:enhanced}
The {\bf enhanced Teichm\"uller space} of $\Sigma$, $\cTh_{g,n}$, 
is the space of $n+1$-uples $(h, \epsilon_1, \cdots, \epsilon_n)$, where $h$
is a hyperbolic metric with geodesic boundary on $\Sigma$ and, for each $k\in
\{ 1, \cdots, n\}$, $\epsilon_k$ is:
\begin{itemize}
\item $0$ if the corresponding boundary component of $\Sigma$ corresponds to a
  cusp of $h$,
\item either $+$ or $-$ if the corresponding boundary component of $\Sigma$
  corresponds to a geodesic boundary component of $h$.
\end{itemize}
\end{defi}

Fock showed in particular that shear coordinates on a surface with
some boundary components provide a natural parametrization of this
enhanced Teichm\"uller space.

\begin{figure}
\begin{center}
\epsfig{figure=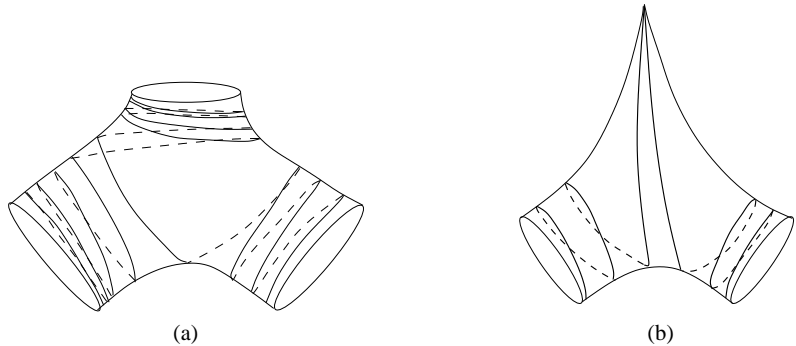,height=1.8in}
\end{center}
\caption{A boundary component can degenerate into a puncture as a result of 
an earthquake. Figure (a) shows the
surface before and (b) after an earthquake.}
\label{fig:pants}
\end{figure}

Note that the boundary of $\cT_{g,n}$ has a stratified structure, with strata
corresponding to subsets of the set of boundary components which are
``pinched'' to obtain cusps, as shown in Figure \ref{fig:pants}. 
Heuristically, $\cTh_{g,n}$ is obtained by ``reflecting'' $\cT_{g,n}$ along
the codimension $1$ strata of its boundary, and $\cTh_{g,n}$
contains an open dense subset which is a $2^n$-fold cover of the interior of
$\cT_{g,n}$. There is also a natural embedding of $\cT_{g,n}$ in $\cTh_{g,n}$,
obtained by taking all $\epsilon_i$ equal to $+$ in the definition above.

It is possible to define in a rather natural -- but perhaps not obvious -- 
way the element of $\cTh_{g,n}$ obtained by an earthquake along a measured
geodesic lamination, i.e., a map $E_r:\cML_{g,n}\times \cTh_{g,n}\rightarrow 
\cTh_{g,n}$. This map has the key properties that should be required of it:
\begin{itemize}
\item its restriction to $\cT_{g,n}$ (considered as a subset of $\cT_{g,n}$), 
followed by the projection from $\cTh_{g,n}$ to $\cT_{g,n}$, is the right
earthquake map $E_r:\cML_{g,n}\times \cT_{g,n}\rightarrow\cT_{g,n}$ defined 
above,
\item it is continuous,
\item for any $\lambda\in \cML_{g,n}$, any $h\in \cTh_{g,n}$ and any
$t,t'\in \R_{>0}$, 
$$ (E_r(t\lambda)\circ E_r(t'\lambda))(h) = E_r((t+t')\lambda)(h)~. $$ 
\end{itemize}

Theorem \ref{tm:earthquake} can then be reformulated in a simpler way in terms of
$\cTh_{g,n}$. 

\begin{teo} \label{tm:earthquake2}
For any $h,h'\in \cTh_{g,n}$, there exists a unique $\lambda\in \cML_{g,n}$ 
such that $h'=E_r(\lambda)(h)$.
\end{teo}

It is shown in section 9 how 
Theorem \ref{tm:earthquake2} follows from Theorem
\ref{tm:earthquake}. Note that some care is needed there to give
the proper definitions and prove the result. 

\subsubsection*{The Mess proof of the Earthquake Theorem.}

G. Mess \cite{mess} discovered some striking similarities between
quasifuchsian hyperbolic 3-manifolds and the so-called GHMC (for ``globally
hyperbolic compact maximal'') AdS (for ``Anti-de Sitter'') 3-dimensional
manifolds. As a consequence he found a direct and very geometric proof 
of the Earthquake Theorem.

The 3-dimensional AdS space, $AdS_3$, can be defined as a quadric in $\R^4$
endowed with a symmetric bilinear form of signature $(2,2)$, with the induced
metric:
$$ AdS_3=\{ x\in \R^{2,2}~|~\langle x,x\rangle=-1 \}~. $$
It a complete Lorentz space of constant curvature $-1$, analog in certain ways
to the hyperbolic 3-space. Defined in this way, $AdS_3$ is however not simply
connected, its fundamental group is $\Z$. Its totally geodesic planes are
isometric to $H^2$, while its time-like geodesics are closed of length $2\pi$.

An AdS manifold is a manifold endowed with a Lorentz metric locally isometric
to the metric on $AdS_3$. Recall that a {\it Cauchy surface} in a Lorentz manifold
is a surface which intersects each inextendible time-like geodesic exactly once,
see e.g. \cite{O}. We are particularly interested here in 
{\it globally hyperbolic maximal compact (GHMC)} AdS 3-manifolds:
those AdS 3-manifolds which contain a closed, space-like
Cauchy surface, and which are maximal under this conditions (any isometric
embedding into an AdS manifold containing a closed Cauchy surface is an
isometry). GHMC AdS manifolds display some striking similarities with
quasifuchsian hyperbolic 3-manifolds.

Mess discovered in particular that the space of
GHMC AdS manifolds which are topologically
$\Sigma \times \mr$ (where $\Sigma$ is a closed surface of genus at
least $2$) is parametrized by the product of two copies of the Teichm\"uller
space of $\Sigma$, $\cT_\Sigma$. This is strongly reminiscent of the Bers
double uniformization theorem \cite{bers}. However it 
does not involve a conformal
structure at infinity, but rather the ``left'' and ``right'' hyperbolic
metrics, $h_l$ and $h_r$, associated to such an AdS 3-manifold (the
definitions can be found in Section 2).

Moreover those GHMC AdS manifolds have a ``convex core'', and the boundary of
this convex core has two connected components, each with an induced hyperbolic
metric (which we call $\mu_+$ and $\mu_-$) and a measured bending lamination
(called $\lambda_+$ and $\lambda_-$ here). The left hyperbolic metric $h_l$ is
obtained from the induced metric on the upper boundary component of the convex
core, $\mu_+$, by the action of the left earthquake relative to $\lambda_+$
(rather than by a grafting along $\lambda_+$, as in the quasifuchsian
context). This leads to the following diagram, where $E_l(\lambda)$ (resp.
$E_r(\lambda)$ is the left (resp. right) earthquake relative to the measured
lamination $\lambda$.

\begin{diagram}
& & \mu_+ & & \\
& \ldTo^{E_l(\lambda_+)} & & \rdTo^{E_r(\lambda_+)} & \\
h_l & & & & h_r \\
& \luTo_{E_r(\lambda_-)} & & \ruTo_{E_l(\lambda_-)} & \\
& & \mu_- & & 
\end{diagram}

It follows that $h_l=E_l(2\lambda_+)(h_r)=E_r(2\lambda_-)(h_r)$. Since any
couple $(h_l, h_r)$ can be obtained as the left and right hyperbolic metrics
of exactly one GHMC AdS manifold, a simple proof of the Earthquake Theorem
follows. 

This line of ideas can be extended to obtain an ``Earthquake Theorem'' for 
hyperbolic metrics with cone singularities, of fixed angle in $(0,\pi)$, on
closed surfaces, see \cite{cone}. The GHMC AdS manifolds considered by Mess
are then replaced by similar manifolds with ``particles'', i.e., cone
singularities along time-like geodesic segments.

\subsubsection*{Multi-black holes.}

There is a class of 3-dimensional AdS manifolds analogous to GHMC
manifolds, which is obtained by replacing the closed Cauchy surface by a non-compact
one. These manifolds were first defined in the physics literature
\cite{Multi,Brill} and are called ``multi-black holes'' (called MBH here). 
A mathematical description can be found in~\cite{barbot-1,barbot-2}. 
The simplest example is obtained from a complete hyperbolic metric $h$ on a
compact surface $S$ of genus $g$ with $n$ disks removed (with each end
of infinite area) by a warped product construction:
$$ M =(S\times (-\pi/2,\pi/2), -dt^2+\cos(t)^2 h)~. $$
More general MBH metrics are obtained by deforming those examples, losing the
symmetry $t\mapsto -t$.
 
It is in particular proved in~\cite{barbot-1,barbot-2}  that, 
given a compact surface with 
boundary $S$, the space of MBHs which are topologically the product of $S$ by
an  interval is parameterized by the product of two copies of the Teichm\"uller
space of hyperbolic metrics with geodesic boundary on $S$, as was proved by
Mess for closed surfaces \cite{mess}.

\subsubsection*{The geometry of multi-black holes and the idea of the proof.}

Let $M$ be an MBH, with fundamental group $\pi_1(\Sigma)$. The main idea of
the proof of Theorem \ref{tm:earthquake} is to consider a special class of
convex pleated surfaces in a MBH. It was proved in \cite{benedetti-bonsante} 
that given a MBH $M$ with right and left holonomies $h_l$ and $h_r$, 
there is a one-to-one correspondence between
\begin{itemize}
\item space-like, convex, pleated, inextendible surfaces in $M$ (in general not
complete, but with geodesic boundary),
\item earthquakes between pairs of hyperbolic surfaces with convex boundary
(of finite or infinite area, possibly with vertices at infinity) with
left and right holonomies equal to $h_l$ and $h_r$.
\end{itemize}
One key technical result here is that, given $M$, there is a finite number
of convex pleated surfaces for which each boundary component is either
a closed geodesic or a cusp. Those surfaces have a simple characterization
in terms of the quotient of the boundary components of the convex hull of some natural
curves complementing the limit set of $M$ in a ``boundary at infinity'' of
$AdS_3$ (see the first paragraph of section 3), as shown in Proposition \ref{pr:51}.

In a previous version of this paper, multi-black holes played a key role in 
the proof of the main result. Here however this proof has been rewritten to
be readable to readers with no previous knowledge of multi-black holes.
Some elements of the geometry of with multi-black holes, and the relation
with the main theorem here, are explained in section 10.

\subsubsection*{A description in terms of measured laminations.}

A by-product of the arguments used for the proof of Theorem
\ref{tm:earthquake} is another description of the space of MBHs of given
topology, based on pleated surfaces or, in other terms, on hyperbolic metrics
and measured laminations on compact surfaces with boundary. 
This is explained in more details in the physics introduction of a previous
version of this text, see \cite{mbh-v4}. We do not dwell on this point here.

\subsection*{Acknowledgements.}

We would like to thank Thierry Barbot and 
Francis Bonahon for some useful conversations and comments.

\section{Earthquakes on $\Tt_{g,n}$}

\subsection{The Teichm\"uller space $\cT_{g,n}$.}

A hyperbolic metric $\eta$ on $\Sigma$ is said  to be \emph{admissible} if:
\begin{enumerate}
\item It has a finite area.
\item Its completion has a geodesic boundary.
\item Each geodesic boundary component is a closed curve.
\end{enumerate}

We denote by $\Sigma_\eta$ the hyperbolic surface $(\Sigma, \eta)$, 
and by $\overline\Sigma_\eta$ the completion of
$\Sigma_\eta$. Notice that the topological type of $\overline\Sigma_\eta$
depends on $\eta$. A neighbourhood of a puncture can look like either a
cusp or a neighborhood of a boundary component.

The Teichm\"uller space $\Tt_{g,n}$ for $\Sigma$ is the space of 
admissible hyperbolic metrics up to the action
of diffeomorphisms isotopic to the identity. For $\chi(\Sigma)<0$
this space is non-empty.

Given an admissible metric on $\Sigma$, its holonomy is a faithful
(i.e. injective) and 
discrete representation
\[
    h:\pi_1(\Sigma)\rightarrow PSL_2(\mr)\,.
\]
The surface $\overline{\Sigma}$ is the convex core $\Kk$ 
of the quotient of 
$\mathbb H^2$ (hyperbolic plane) by
the action of $\Gamma:=h(\pi_1(\Sigma))$. One can easily check that the
following  statement holds\\

\emph{For each $\gamma\in\pi_1(\Sigma)$ parallel to a puncture, either
  $h(\gamma)$ is parabolic  or its axis is a boundary curve of $\Kk$.}  (*)\\

A  faithful and discrete representation
$h:\pi_1(\Sigma)\rightarrow PSL_2(\mr)$ satisfying (*) is
called~\emph{admissible}. Thus, the holonomy of an
admissible metric is an admissible representation. 
Conversely, the quotient of the convex core
of an admissible representation is  a finite area
hyperbolic surface homeomorphic to $\Sigma$.
Thus, the space $\Tt_{g,n}$ can be identified with the space of admissible
representations of $\pi_1(S)$ into $PSL_2(\mr)$, up to conjugacy. 

Since the fundamental group of $\Sigma$ is a free group on $2g+n-1$ generators
it follows that the space of representations of $\pi_1(\Sigma)$ into
$PSL_2(\mr)$ is $PSL_2(\mr)^{2g+n-1}$. Taking into account the fact that
conjugate 
representations lead to the same metrics we see that $\dim\Tt_{g,n}=6g-6+3n$.
The Teichm\"uller space $\Tt_{g,n}$ is a closed subset of
this space with interior  corresponding exactly to the metrics without
cusps. The boundary of $\Tt_{g,n}$ corresponds to structures with some cusps. 
 
\subsection{Measured geodesic laminations on a hyperbolic surface with geodesic
  boundary}

Let us fix an admissible metric $\eta\in\Tt_{g,n}$ with holonomy
$h:\pi_1(\Sigma)\rightarrow PSL_2(\mr)$.

A  geodesic lamination on $\Sigma_\eta$ is a closed subset
$L$ foliated by complete geodesics. A \emph{leaf} of $L$ is a geodesic of the
foliation, whereas a \emph{stratum} is either a leaf or a connected component
of $\Sigma_\eta\setminus L$.

Since the area of $\Sigma_\eta$ is finite, the structure of
$L$ can be proved to be similar to the structure of a geodesic lamination 
on a closed surface. In particular:
\begin{itemize}
\item The Lebesgue measure of $L$ is $0$.
\item There exists a unique partition of $L$ in complete geodesics (that is,
  the support $L$ is sufficient to encode the lamination).
\item $\Sigma_\eta\setminus L$ contains finitely many connected
  components. Each of them is 
  isometric to (the interior of) a finite area hyperbolic surface with
  geodesic boundary. 
\end{itemize}
A leaf of $L$ is a \emph{boundary curve} if it is the boundary of some
component of $\Sigma_\eta\setminus L$.
\begin{itemize}
\item Boundary curves are finitely many. Moreover they are dense in $L$.
\end{itemize}

The following lemma describes the behaviour of a geodesic lamination near a
puncture. 

\begin{lem} \label{lm:31}
For each boundary component $c$ there exists an $\eps$-neighbourhood $U$ 
such that every leaf intersecting $U$ must spiral around $U$.
Moreover, leaves in $U\cap L$ are locally isolated.

The same result holds for cusps, by exchanging $\eps$-neighbourhoods by
horoballs: for each cusp $c'$ there exists a neighborhood $U$ bounded by a 
horocycle $C$ such that 
every leaf intersecting $C$ does so orthogonally, and leaves in $U\cap C$ are
locally isolated.
\end{lem}

\pf
We prove the first part of the statement. The case with cusp is completely
analogous. 
On the other hand the proof uses the same arguments used in~\cite{Casson}
to describe the behaviour of a geodesic lamination (without measure) on a closed surface 
in a regular neighbourhood of some closed leaf.

Let $\Sigma_\eta=\Hh/h$ where $h$ is
the holonomy representation of $\pi_1(\Sigma)$ and $\Hh$ is the convex core of
$h$.

Let $\tilde L$ be the pre-image of $L$ on $\mh^2$ , $\tilde c$ be a pre-image
of $c$ and $\gamma$ be a generator of the stabilizer of $c$. 
If $d$ is the  length of $c$, we may find $\eps>0$
such that if $\tilde c'$ is a geodesic $\eps$-close to $\tilde c$ 
then the length of
the projection of $\tilde c'$ on $c$ is greater  than $d$. 
Thus if $\tilde c'$ is at positive distance from  $\tilde c$ then
$\gamma\tilde c'$ must intersect $\tilde c'$.

Thus leaves of $L$ intersecting $U_\eps$ have to spiral around $c$.

Now let us prove that leaves in $U\cap L$ are locally isolated. By taking a
smaller $\eps$, we may suppose that $U_\eps$ projects on a regular
neighbourhood of $c$.
Take a leaf spiraling around $c$, say $l$, and denote by $\tilde l$ a lifting
of $l$ on $\mh^2$ intersecting $U_\eps$. 
Suppose that between $l$ and $\gamma l$ there are
infinitely many leaves intersecting $U_\eps$. 
Thus there are infinitely many boundary leaves. On the
other hand that leaves  between $l$ and 
$\gamma l$ intersecting $U_\eps$ are not permuted by
$\pi_1(\Sigma)$, so we get a contradiction. 
(It follows from this
argument that there are finitely many boundary
leaves in $\Sigma$.)
\cvd

\subsubsection*{Transverse measures.}

The notion of transverse measure can be introduced as in the closed case.
We say that an arc in $c$ is transverse to $L$ if it is transverse to 
the leaves of $L$.

A transverse measure on $L$ is the assignment of a Borel measure $\mu_c$ on
each transverse arc $c$ such that:
\begin{enumerate}
\item The support of $\mu_c$ is $c\cap L$.
\item If $c'\subset c$ then $\mu_{c'}=\mu_c|_{c'}$.
\item If two transverse arcs are homotopic through a family of transverse arcs
   then their total masses are equal.
\end{enumerate}

The simplest example of a geodesic lamination is a simple geodesic $u$. In such a
case 
a measure $\mu_c$ is concentrated on the intersection points of $c$ with
$u$. The mass of each single intersection point is a number independent of $c$
and 
is, by definition, the weight of $u$. Thus transverse measures on $u$ are
encoded by a positive number. 

On closed surfaces, every measured geodesic lamination splits as the disjoint 
union of sub-laminations
\[
   L=S\cup L_1\cup L_2\cup\ldots...\cup L_k
\]
such that the support of $S$ is a finite union of simple geodesics and each
leaf $l\subset L_i$ is dense in $L_i$.

In the case we are concerned with, things are a bit more complicated, since $L$
is not supposed to be compact. On the other hand we have seen that  near a 
puncture $L$ has a simple behaviour.
Notice that a consequence of Lemma~\ref{lm:31} 
is that every geodesic in $L$ that enters a
cusp or spirals around a geodesic boundary is weighted. Thus it cannot have
accumulation points in $\Sigma$. 
It follows that  such leaves are properly embedded in $\Sigma$. So, if some
regular neighbourhoods of the punctures are cut off from $\Sigma$, such leaves 
appear as properly embedded compact arcs.

This remark allows to find
a canonical decomposition of a measured geodesic lamination.

\begin{lem}\label{a:lem}
If $\lambda$ is a measured geodesic lamination on $(\Sigma,\mu)$, then it
splits as the union of sub-laminations
\[
   L=B\cup S\cup L_1\cup L_2\cup\ldots\cup L_k
\]
such that $B$ is the union of leaves that do not have
compact closure in $\Sigma$, 
$S$ is a union of closed geodesics. $L_i$ is compact and every leaf
$l$ of $L_i$ is dense in $L_i$.
\end{lem}

\pf
Define first $B$ as the union of the geodesics in the support of $\lambda$
that enter any neighbourhood of the boundary. Their behaviour near the boundary
is described by Lemma \ref{lm:31}. 
Let $\lambda'$ be the measured lamination
obtained by removing from $\lambda$ the measure supported on $B$.

We now consider the surface $(\Sigma',\mu')$ obtained by gluing 
two copies of $(\Sigma, \mu)$ along their boundary, by identifying
corresponding points of the boundary on the two copies. Since the support of
$\lambda'$ does not enter some neighbourhood of $\dr\Sigma$, $\lambda'$ lifts
to a measured geodesic lamination on $(\Sigma',\mu')$. Applying the known
decomposition result for closed surfaces to $\lambda'$ on $\Sigma'$ shows that
its support can be written as $S\cup L_1\cup \cdots \cup L_k$, and the
result for $\lambda$ follows. 
\cvd

Measured geodesic laminations with compact support are well understood.
To get a complete description of a general measured geodesic lamination,
 we should  describe complete embedded geodesics of
$\Sigma$ that escape from compact sets.

We have seen that every leaf $l$ in $B$ produces a properly embedded arc in the
complement of some regular neighbourhood of the puncture. Notice that the
homotopy class of this  arc does not depend on the regular
neighbourhood. With a slight abuse of language we say that $l$ represents
such a class.

We could expect that $l$ is determined by its homotopy class. This is not
completely true. In fact the homotopy class does not ``see'' in which way $l$
winds around the boundary of $\Sigma_\eta$.

\begin{lem}\label{b:lem}
Suppose that a positive way of spiraling around each boundary component of
$\overline\Sigma_\eta$ is fixed. Then in each homotopy 
class of properly embedded arcs joining
two puncture of $\Sigma$, there exists a unique geodesic representative that
spirals in the positive way.
\end{lem}

\pf
Let $c_1$ and $c_2$ be two punctures of $\Sigma$, and let $h$ be a homotopy
class of properly embedded arc joining them. $c_1$ and $c_2$ correspond to
geodesic boundary components of $\Sigmab_\eta$, which we still call $c_1$ and
$c_2$. Let $c'_1$ be a lift of $c_1$ as a connected component of the
(geodesic) boundary of the universal cover of $\Sigmab_\eta$, and similarly
let $c'_2$ be a lift of $c_2$ as a connected component of $\dr\Sigmab_\eta$,
chosen so that there is a lift $h'$ of $h$ as a path connecting $c'_1$ to
$c'_2$. 

Any realization of $h$ as a geodesic spiraling around $c_1$ and $c_2$ has to
lift to the universal cover of $\Sigmab_\eta$ as a geodesic which is
asymptotic to $c'_1$ and $c'_2$. There are four such geodesics, depending on
the choice of one of the two ends of $c'_1$ and one of the two ends of $c'_2$.
But only one of those choices corresponds to the positive spiraling direction,
so there is only one geodesic realization of $h'$.
\cvd 

Given an admissible metric $\eta$, denote by $\Mm\Ll_{g,n}(\eta)$ the set of
measured geodesic laminations on the surface $\Sigma_\eta$. From Lemmas~\ref{a:lem}
and \ref{b:lem} it follows that if $\Sigma_\eta$ and $\Sigma_{\eta'}$ have no
cusp, then there is a natural bijection
\begin{equation}\label{a:eq}
    \Mm\Ll_{g,n}(\eta)\rightarrow\Mm\Ll_{g,n}(\eta')\,.
\end{equation}
Actually given a measured geodesic lamination $\lambda$ on
$\Mm\Ll_{g,n}(\eta)$ it is the union of a compact sub-lamination $\lambda_c$
and a sub-lamination $\lambda_b$ of leaves spiralling along some bounary
components. 
Now, there is a compact measured geodesic lamination $\lambda'_c$ 
in $\Mm\Ll_{g,n}(\eta')$ obtained by ``straightening'' leaves of $\lambda_c$ with
respect to $\eta'$ (it is possible for instance to consider $\Sigma$ as 
included in its double and apply the analogous result for laminations in a
closed surface). Moreover by Lemma \ref{b:lem} we can also straighten the
lamination $\lambda_b$ with respect to $\eta$, and the union of
$\lambda'_c\cup\lambda'_b$ corresponds to $\lambda$ via identification
(\ref{a:eq}).

When $\eta'$ is supposed to have some cusps, the map~(\ref{a:eq}) can be defined
in the same way, but it is no longer $1$-to-$1$. The reason is that if we change the
orientation of spiralling of leaves along a geodesic boundary of $\eta$ that is a
cusp of $\eta'$, the corresponding lamination of $\eta'$ does not change at
all.

In this work we will denote by $\Mm\Ll_{g,n}$ the set of measured geodesic
laminations of a hyperbolic surface with geodesic boundary (without
cusps). From the above discussion this set is well-defined and for every
admissible metric $\eta$ we have a surjective map
\[
 \Mm\Ll_{g,n}\rightarrow\Mm\Ll_{g,n}(\eta)\,.
\]

\subsection{The mass of boundary component.}

Given a measured geodesic lamination $\lambda$ on $\Sigma_\eta$, 
the \emph{mass} of a puncture with respect to $\lambda$ is a positive number
$m_\lambda(c)$ that measures how much the measured 
lamination is concentrated in a neighbourhood of $c$. 

We will give the construction of $m_\lambda(c)$, when $c$ corresponds to a
geodesic boundary component of $\overline\Sigma_\eta$. 

Fix a regular neighborhood $U_\eps$ of $c$ such that
every leaf intersecting $U_\eps$ spirals around $c$.
For every $x\in U_\eps$ consider the geodesic loop $c_x$ with vertex at $x$ parallel
to $c$. We claim that the total mass of such a loop does not depend on $x$.

Let $\Hh$ be the convex core of the holonomy $h$ of $\Sigma_\eta$.
Choose a lifting
of $c$, say $\tilde c\subset\partial\Hh$ and let $\gamma$ be the generator of the stabilizer of
$\tilde c$ in $\pi_1(\Sigma)$. 
If $\tilde x$ is a lifting of $x$ then the loop $c_x$ lifts to the
segment $[x,h(\gamma) x]$. Since geodesics spiraling around $c$ lift to geodesics
asymptotic to $\tilde c$, it follows that $c_x$ intersects every such geodesic
once.
Since the total mass of $c_x$ depends only on the number of intersection
points of $c_x$ with each leaf, it does not depend on $x$.

The same construction works when $c$ corresponds to a cusp.

Notice that $m_\lambda(c)=0$ if and only if there exists a neighbourhood of $c$
avoiding $L$.

When $c$ corresponds to a geodesic boundary, the total mass of $c$ does not
give information about the orientation of spiraling of leaves around $c$.
If we choose for each boundary component a positive way of spiraling, then we
can define a \emph{signed} mass of $\overline m(c)$ in the following way:
\begin{itemize}
\item $|\overline m(c)|=m(c)$;
\item $\overline m(c)>0$ if and only if it spirals in the positive way around $c$.
\end{itemize}
(The second requirement makes sense because two leaves near $c$ have to spiral
in the same way.)

Let us stress that the signed mass of $c$ can be defined only for punctures
corresponding to geodesic boundary components, and it is well defined up to the choice of a
positive way of spiraling.

\subsection{Geodesic laminations on a pair of pants.}

Here we give an explicit description of the measured geodesic laminations on a
hyperbolic pair of pants in terms of the signed masses. This case is relevant
to what in the physics literature is known as the 3 asymptotic region black hole
(see \cite{barbot-1,barbot-2}).

\begin{prop} \label{pr:pants-lamin}
Fix a hyperbolic pair of pants $P$, and for each boundary component choose a positive way 
of spiraling.
Then the function that associates to every measured geodesic lamination on 
$P$ the signed masses of the boundary components of $P$ is bijective.
\end{prop}

\pf
Denote by $c_1,c_2,c_3$ both the punctures of $\Sigma_{0,3}$ and the
corresponding boundary curves on $P$.

Since simple closed curves in $\Sigma_{0,3}$ are boundary parallel, geodesic
laminations do not contain a compact part. 
Moreover notice that there are $6$ properly embedded arcs up to
homotopy. Each of them is determined by its end-points. There are three arcs
connecting different punctures and three arcs connecting the same puncture.

Thus there exist exactly $4$ maximal systems of disjoint properly
embedded arcs in $\Sigma_{0,3}$.  Namely, $L_0$ is the union of arcs
connecting different components whereas $L_i$ (for $i=1,2,3$) is the union of arcs
with endpoint at $c_i$.

Given three positive numbers $m_1, m_2, m_3$, an explicit computation  shows 
that only one of $L_i$ can be equipped with a system of weights which give
masses equal to $m_i$. The system of weights is uniquely determined as well.

In particular the measures on $L_0$ correspond to $m_1,m_2,m_3$
satisfying triangular inequalities, whereas measures on $L_i$ correspond to the case
$m_i\geq m_j+m_k$.
\cvd

\section{Earthquakes}

In this section we recall the definition of earthquakes on hyperbolic surfaces,
in a way which is adapted to hyperbolic surfaces with geodesic boundary, and
show how the definition can be extended to this setting.

\subsection{Earthquakes on convex subsets of $\mh^2$ with geodesic boundary}

Let $\mathcal H$ be an \emph{open} convex set with geodesic boundary
in $\mathbb H^2$ and $L$ be a geodesic lamination of $\mathcal H$.
By definition, a stratum of $L$ is either a leaf of $L$ or 
a component of $\mathcal H\setminus L$.
A right earthquake on $\mathcal H$ with fault locus $L$
is a (possibly discontinuous) map
\[
  E:\mathcal H\rightarrow\mh^2
\]
with the property that 
\begin{itemize}
\item for every stratum $F$, there is an isometry $A(F)\in PSL_2(\mathbb R)$
such that $E|_F=A(F)|_F$,
\item given two strata $F$ and $F'$ the comparison map $A(F)^{-1}\circ A(F')$ is
a hyperbolic transformation whose axis weakly separates $F$ from $F'$ and translates
$F'$ to the right as seen from $F$.
\end{itemize}

Given an earthquake on $\mathcal H$ with fault locus $L$, we can equip $L$ with 
a transverse measure that encodes the amount of shearing.
More precisely given a path $c:[0,1]\rightarrow\mathcal H$ transverse to $L$ and given a partition
$I=(0=t_0<t_1<\cdots<t_k=1)$ we consider the number $\mu(c;I)$ that is the sum
of the translation lengths of the comparison maps $A(F(t_{i+1})A(F(t_i))^{-1}$ where
$F(t)$ is the stratum through $c(t)$.

By a standard fact of hyperbolic geometry on the composition of hyperbolic 
transformations with disjoint axes, if $I'$ is finer than $I$ then
$\mu(c;I')\leq\mu(c;I)$. Thus we can define
\[
   \mu(c)=\inf_{I}\mu(c;I)=\lim_{|I|\rightarrow 0}\mu(c;I)~,
\]
and $\mu$ defines a transverse measure on $L$.

Thurston showed that the measured lamination $\lambda=(L,\mu)$ determines the earthquake $E$ 
(\cite{thurston-earthquakes}).

\begin{prop}\label{earth-lam:prop}
Given a measured geodesic lamination $\lambda$ on $\mathcal H$, there is a unique
earthquake (up to post-composition with isometries of $\mathbb H^2$)
with shearing lamination $\lambda$.
\end{prop}

Contrary to the case discussed in \cite{thurston-earthquakes} 
where earthquakes are bijective maps from $\mh^2$ to itself,
in our setting the image of the earthquake does not need to be the whole
$\mh^2$. This is the reason why Proposition \ref{earth-lam:prop}
holds in our setting whereas it was not true in \cite{thurston-earthquakes}. 

On the other hand it is not difficult to prove that for every earthquake
$E:\mathcal H\rightarrow\mh^2$ the image $E(\mathcal H)$  
is a convex set with geodesic boundary because it is a connected union of
geodesics and ideal hyperbolic polygons (see Lemma \ref{conn:lm}).

\subsection{Earthquakes on $\Tt_{g,n}$}

Given an admissible hyperbolic metric $\eta$ on $\Sigma$, the
left and right earthquakes along a measured geodesic lamination $\lambda$ can
be defined like in the compact case.
 
When the lamination is locally finite they can be described in a very
simple way.  The right earthquake along $\lambda$ is obtained by
shearing each component of $\Sigma\setminus\lambda$ to the right of the
adjacent component by a factor equal to the mass of the boundary.

For the general case it is convenient to construct an equivariant  earthquake on
the universal covering. 

The universal covering of $\Sigma_\eta$, say $\mathcal H$, is an open convex
subset with geodesic boundary in $\mathbb H^2$.  More precisely
$\mathcal H$ is the convex hull of the limit set of the
holonomy $h$ of $\eta$.

The lifting of $\lambda$ is a $h$-invariant measured geodesic 
lamination $\tilde\lambda$.
Consider the right earthquake along $\tilde\lambda$, say
\[
   E:\mathcal H\rightarrow\mathbb H^2.
\]
By the invariance of $\tilde\lambda$, it turns out that
$E\circ h(\gamma)$ is still an earthquake with
shearing lamination $\lambda$.

By Proposition~\ref{earth-lam:prop}, for every
$\gamma\in\pi_1(\Sigma)$ there is an element $h'(\gamma)\in
PSL_2(\mathbb R)$ such that
\[
    E\circ h(\gamma)=h'(\gamma)\circ E~.
\]

\begin{prop}
The representation $h'$ is faithful and discrete.  
The quotient $\mathbb H^2/h'$ is homeomorphic to $\Sigma$.
The map $E$ induces to the quotient a piece-wise isometry
\[
     E^r_\lambda:\Sigma_{\eta}\rightarrow E(\mathcal H)/h'
\]
The surface $E(\mathcal H)/h'$ concides with the convex core of
$\mathbb H^2/h'$ (it is in particular an admissible surface).
\end{prop}

\begin{proof}
First notice that $h'$ is discrete. Indeed let $p$ be some point contained in the interior of some
2-dimensional stratum $F$ of $\tilde\lambda$. Now the $h'$-orbit of $E(p)$ 
accumulates at $E(p)$ if and only if the $h'$-orbit of $p$ accumulates at $E(p)$.
This show that the orbit of $E(p)$ is discrete. Thus $h'$ is a discrete representation.
Since the earthquake map is injective, it turns out that $h'$ is faithful.

To prove that $\mh^2/h'\cong\Sigma$, notice that $h$ and $h'$ are connected by a path
of faithful and discrete representations. 
Namely let $h_t$ be the representation corresponding to the earthquake along $t\lambda$.

To conclude the proof we have to check that $E(\mathcal H)$ is the
convex hull of the limit set of $h'$.  
Let $\tilde U$ be the lifting on $\mathcal H$ of a regular neighbourhood of punctures in $\Sigma_\eta$.
A simple argument shows $\mathcal H\setminus\tilde U$
is sent by $E$ to a subset with compact quotient.

Thus, it is
sufficient to show that there is a constant $M$, such that 
for any point $p$ close to a puncture $x$
there exists a loop centered at $p$, parallel to the puncture, whose
length is bounded by  $M$.

Take the geodesic loop $\gamma$ of $\Sigma_\eta$ centered at $p$ and parallel to $x$.
Notice that $\gamma$ meets only a finite number of leaves of $L$.

The image of $\gamma$ via $E_\lambda^r$ is a union of geodesic arcs $\gamma_i$  
whose end-points $x_i,y_i$ lie on $E_\lambda^r(L)$. The piece-wise
geodesic loop
\[
   \hat\gamma=\gamma_0 *[y_0,x_1]*\gamma_1*[y_1,x_2]*\ldots*\gamma_N
\]
is parallel to $x$. Notice that the sum of the lengths of $\gamma_i$ is equal
to the length of $\gamma$, whereas the length of the segment $[y_i,x_{i+1}]$
is equal to the mass of the corresponding leaf. Thus the length of
$\hat\gamma$ is equal to the sum of the length and the mass of $\gamma$. 
\end{proof}

We say that $E(\mathcal H)/h'$ is obtained by a right earthquake
of $\Sigma_\eta$ along $\lambda$ and we denote it by $E^r_\lambda(\Sigma_\eta)$.

We have seen that a lamination on $\Sigma_\eta$ is the disjoint union
of a compact part, say $\lambda_c$, and a finite union of leaves that
spirals around boundary components or enter cusps, say $\lambda_b$.
The earthquake along $\lambda$ can be regarded as the composition of
the earthquake along $\lambda_c$ and the earthquake along $\lambda_b$:
more precisely we have to compose the earthquake along $\lambda_c$
with the earthquake along $\hat\lambda_b$ that is the image of
$\lambda_b$ in $E^r_{\lambda_c}(\Sigma_\eta)$.

The earthquake along $\lambda_c$ can be easily understood: we
approximate $\lambda_c$ by weighted multicurves. Then the earthquake
along $\lambda_c$ is the limit of the fractional Dehn twists along
these weighted multicurves.
Notice that the earthquake along $\lambda_c$ does not change the
length of any boundary component.

The earthquake along $\lambda_b$ can be described in the following
way.  We cut the surface (only the interior of $\Sigma_{\eta}$) along
the leaves of $\lambda_b$ and we get a surface $\hat\Sigma$ with
geodesic boundary.  Since $\lambda_b$ is locally finite in
$\Sigma_\eta$, every leaf of $\lambda_b$ corresponds to exactly
two boundary components of $\hat\Sigma$.  Then we glue back the
boundary components corresponding to the same leaf $l$, composing the
original glueing with a right translation of factor equal to the
weight of $l$.

Opposite to the previous case, the earthquake along
$\lambda_b$ changes the length of the boundary components (and may transforms cusps
in geodesic boundary components). In the next section we determine
the length of a boundary component after the earthquake.

\subsection{Boundary length and spiraling orientation after an earthquake.}

The mass of a boundary component $c$ for a measured lamination $\lambda$ is in direct
relation with the variation of the length of $c$
under an earthquake along $\lambda$, and also with the way $\lambda$
spirals on $c$. Indeed the image $\lambda'$ of $\lambda$ by the right 
earthquake $E^r_\lambda$
is well-defined, but it might spiral on $c$ differently from $\lambda$. 

Let us choose an explicit way of spiralling around each boundary curve in the
following way. An orientation is induced by $P$ on its boundary. 
If $l$ spirals around some $c_i$ then an orientation is induced on $l$ by the
orientation on $c_i$. Namely $l_i$ is oriented in such a way that the nearest-point retraction on $c_i$
(that is well-defined in a neighbourhood of $c_i$) is orientation preserving. 
Notice that if $l$ spirals around $c_i$ and $c_j$
the orientations induced on $l$ may disagree. 

Then we say that $l$ spirals in a positive way around $c_i$ if it goes closer and
closer to $l$. We call it the \emph{standard} spiraling orientation, and we
will refer to it througout this paper.

\begin{prop} \label{pr:length}
Let $a$ be the length of $c$ in $\Sigma_\eta$,and let $a'$ be the length of the
 corresponding boundary component after a left earthquake along $\lambda$, 
in $E^r_\lambda(\Sigma_\eta)$.
\begin{enumerate}
\item $a'=a+m$ if $\lambda$ spirals around $c$ in the positive way, $a'=|a-m|$
if $\lambda$ spirals around $c$ in the negative way.
\item If $\lambda$ spirals in the positive direction, so does
  $\lambda'$. If $\lambda$ spirals in the negative direction, then
  $\lambda$ spirals in the negative direction if $m<a$, in the
  positive direction if $m>a$.
\item $E^r_\lambda(\Sigma_\eta)$ has a cusp at the boundary component
  corresponding to $c$ if and only if $\lambda$ spirals in the
  negative direction at $c$ and its mass $m$ is equal to the length
  $a$ of $c$ in $\eta$.
\end{enumerate}
\end{prop}

\begin{proof}
Let us consider the lifting of $E^r_\lambda$ to the universal covering
\[
    E:\mathcal H\rightarrow\mathcal H'~.
\]
Let $\tilde c$ be a lifting of $c$. We can choose coordinates on
$\mh^2$ such that $\tilde c$ is the geodesic from $0$ to $\infty$ and
$\mathcal H$ is contained in the region $\{(x,y)|x<0,y>0\}$.

Suppose that $\lambda$ spirals in the positive way around $c$.  This
means that there is an $\epsilon$-neighbourhood $U$ of $\tilde c$ such
that every leaf intersecting $U$ goes to $\infty$.

Let $\gamma\in\pi_1(\Sigma)$ be a positive representative of the
peripheral loop around $c$ .  We have that $h(\gamma)$ can be written 
in a suitable basis as $\begin{pmatrix}
e^a & 0\\0 & e^{-a}\end{pmatrix}$.  Fix a stratum $F$ intersecting
$U$, then the comparison isometry between $F$ and $h(\gamma)(F)$ is
the composition of hyperbolic translations with attractive fixed point
equal to $\infty$. It is not difficult to see that the comparison
isometry is represented by the $SL_2(\mr)$ matrix $\begin{pmatrix} e^m
  & *\\0 &e^{-m}\end{pmatrix}$.  Since $h'(\gamma)$ is the composition
of the comparison isometry with $h(\gamma)$.  This shows that the
translation length of $h'(\gamma)$ is $a+m$.

Moreover notice that $\infty$ is the attractive fixed point of $h'(\gamma)$
and that $E^r(l)$ ends at $\infty$ for every leaf $l$ that ends at
$\infty$. This show that the image lamination spirals in positive way
around $c$.

The other cases can be obtained by similar computations.
\end{proof}


The computation of the proof of Proposition \ref{pr:length}
can also be found in \cite{thurston-notes} in the special case of a
pair of pants, and in \cite{bonahon-liu} in the slightly different setting
of shear coordinates.  The same proposition also holds -- with positive 
and negative orientations reversed -- for a left earthquake.

\subsubsection*{Earthquakes on a pair of pants.}

One could wonder whether the analog of Theorem~\ref{classical:teo} 
holds also for $\Tt_{g,n}$, that is, given $F,F'\in\Tt_{g,n}$ there exists a unique
$\lambda\in\Mm\Ll_{g,n}$ such that the left earthquake along $\lambda$
transforms $F$ into $F'$.
A classical example due to Thurston shows that this is not the case on
a hyperbolic pair of pants.
In this section we will focus on that example. Since explicit computations
are possible we get a complete picture about earthquakes. In the next
sections we will see that the same picture, suitably expanded, 
holds for general surfaces.

Let $\Sigma$ be the thrice-punctured sphere and let $c_1,c_2,c_3$ denote the
punctures.
It is well known that a hyperbolic metric with geodesic boundary on $\Sigma$ is
determined by three positive numbers $a_1,a_2,a_3$ corresponding to the
lengths of the three boundary components. Moreover when $a_i\rightarrow 0$,
the corresponding geodesic boundary component degenerates to a cusp.
Thus $\Tt_{0,3}$ is parametrized by a triple of non-negative numbers. Let
$P(a_1, a_2, a_3)$ denote the element of $\Tt_{0,3}$ corresponding to the
triple $(a_1, a_2, a_3)$.


We have seen in Proposition \ref{pr:pants-lamin}
that each measured geodesic lamination on $P$ is determined by
three real numbers (the signed masses with respect to the standard spiralling
orientation). Denote by $\lambda(m_1,m_2,m_3)$ the lamination corresponding to
the triple $m_1,m_2,m_3$.  Then the surface obtained by the right earthquake
along $\lambda(m_1,m_2,m_3)$ on $P(a_1,a_2,a_3)$ is
\[
     P(|a_1+m_1|, |a_2+m_2|, |a_3+m_3|)
\]
whereas the surface obtained by a left earthquake is
\[
    P(|a_1-m_1|, |a_2-m_2|, |a_3-m_3|)\,.
\]
Notice that this formulas makes sense also when some $a_i=0$.
In fact in such a case they depend only on $|m_i|$ (we have previously
remarked that it is not possible to define a signed mass corresponding to a
cusp).

It follows from those formulas that two hyperbolic pairs of pants (without cusps) are related by $8$
earthquakes. In fact for each $i$ we can choose in arbitrary way
the corresponding sign of $m_i$.

Let us focus on some points.
\begin{enumerate}
\item Given two hyperbolic pairs of pants $P_0=P(a_1,a_2,a_3)$ and $P_1=P(b_1, b_2, b_3)$
there exists a unique lamination $\lambda$ such that $E^r_{\lambda}(P_0)=P_1$
and the path $E^r_{t\lambda}(P_0)$ is contained in the interior of
$\Tt_{3,0}$, for $t\in [0,1]$. 
Namely $\lambda=\lambda(b_1-a_1, b_2-a_2, b_3-a_3)$.

\item  Take a measured geodesic lamination $\lambda=\lambda(m_1, m_2, m_3)$ and
suppose $m_1<-a_1$. Consider the earthquake path
\[
      P_t=E^r_{t\lambda}P(a_1, a_2, a_3)\qquad \lambda_t=\Ee_{t\lambda}(\lambda)
\]
for $t\in [0,1]$.  It has a critical value at $t_0=-m_1/a_1$. Let us give a
picture of the behaviour of $P_t$ near $t_0$. For $t=t_0-\eps$
the geodesic boundary $c_1$ is very small and by consequence there is a
``big'' regular neighbourhood $U$ (that is the distance of $\partial
U$ from $c_1$ is big). The geodesic lamination spirals in the
positive direction, but it looks almost unwind.
At time $t_0$ the geodesic boundary has disappeared and we have a cusp. The
geodesic lamination is completely un-winded. As $t$ becomes greater than $t_0$,
$c_1$ turns out to be a geodesic boundary component, but this time the geodesic lamination spirals in the
opposite direction.

\item Let $\cTh_{3,0}$ be the space of admissible hyperbolic structures on
   $\Sigma$ equipped with a positive spiraling orientation on each boundary
   component, that is, the enhanced Teichm\"uller space of the thrice-punctured sphere. 
   Notice that this space could be identified with $\mr^3$. In fact
   each such surface is determined by three non-negative numbers (the
   lengths) and a certain number of ``signs'' corresponding to non-zero
   numbers. $\Tt_{3,0}$ could be regarded as the quotient of
   $\cTh_{3,0}$ by the action of the group $G=(\mz^2)^3$, generated by
   the symmetries along coordinates planes. 

\item Take an element $P_0=P(a_0, a_1, a_2)\in\cTh_{3,0}$ (notice that
     $a_i\in\mr$). 
     Given a measured geodesic lamination $\lambda$ with signed mass (with
     respect to the spiraling orientation of 
     $P_0$) equal to $m_1, m_2, m_3$ we can notice that the signed
     masses with respect to the  \emph{canonical} spiraling orientation are given
     by $$ \hat{m}_i= sign(a_i) m_i.$$
\item Suppose all $a_i\neq 0$ (that is $P_0$ is a pair of pants). 
     Then a spiraling orientation can be pushed forward on the surface
     $E^r_\lambda(P_0)$: namely  a curve $c\subset E^r_\lambda(P_0)$ spirals in
     positive direction around $c_i$ iff so does its pre-image in $P_0$.
     Thus earthquakes ``lift'' to a map 
     \[
       \Ee_\lambda:\cTh_{3,0}\setminus\{\textrm{structures with
         cusp}\}\rightarrow \cTh_{3,0}.
     \] 
     If $b_1,b_2,b_3$ are the real parameters corresponding to
     $\Ee_\lambda(P_0)$ we have
     \[
       \begin{array}{l}
         |b_i|=||a_i|-\hat{m}_i|=|a_i-m_i|\\
         sign(b_i)=sign(a_i)sign(|a_i|-\hat{m}_i)=sign(a_i-m_i)
       \end{array}
     \]
     so we get the simple formula
     \[
          b_i= a_i-m_i~.
     \]
     In particular $\Ee_\lambda$ extends on the whole of $\cTh_{0,3}$. Notice
     that if $P_0$ has a cusp in $c_1$ then the orientation of
     $\Ee_\lambda(P_0)$ in $c_1$ depends only on the sign of $m_i$.

\item $\Ee_\lambda$ is not $G$-equivariant on $\cTh_{3,0}$. On the other hand it
  is uniquely determined by the following conditions:
  \begin{itemize}
  \item if  $\Tt_{3,0}$ is identified with the subset of  $\cTh_{3,0}$ corresponding
  to triples $(a_1,a_2,a_3)$ with $a_k\geq 0$, then $\pi\circ
  \Ee_\lambda=E^r_\lambda$ (where $\pi:\cTh_{3,0}\rightarrow\Tt_{3,0}$ is the
  projection).
\item $\Ee_\lambda$ is a flow, that is $\Ee_{t\lambda}\circ
    \Ee_{t'\lambda}=\Ee_{(t+t')\lambda}$ for $t,t'>0$.
  \end{itemize}
\item On $\cTh_{0,3}$ the earthquake theorem holds. That is there exists a unique
   right earthquake joining two points in $\cTh_{3,0}$.
\end{enumerate}


In the next sections we will see that the same picture holds for any surface
$\Sigma_{g,n}$. The key ingredient to get such a generalization is to relate
earthquakes to bent surfaces in the multi-black holes that are defined
in the next section. The relation between earthquakes
and bent surfaces will be obtained by generalizing the Mess argument in the closed
case. The main difference will be that in a multi-black hole there are (finitely) many bent
surfaces (in contrast in the closed case where there is a unique one).

\section{The geometry of Anti de Sitter space}

We collect in this section, for the reader's convenience, some basic facts on the
geometry of the 3-dimensional AdS space, as can be found in particular in 
\cite{mess,mess-notes}.

\subsubsection*{The AdS space and its conformal boundary}

 Let $\mathbb R^{2,2}$ denote $\mathbb R^4$ equipped with
 the standard bilinear symmetric form, say $\langle\cdot,\cdot\rangle$, of
 signature $(2,2)$.

Let us consider the set of negative unit vectors:
\[
  X:=\{x\in\mathbb R^{2,2}~|~\langle x,x\rangle=-1\}\,.
\]
Since the tangent plane at $x$ of $X$ is the linear plane orthogonal to
$x$ with respect to $\langle \cdot,\cdot\rangle$, the induced
symmetric form on $X$ has Lorentzian signature.

The projection
\[
   \pi:X\rightarrow\mpi^3
\]
is a $2:1$ covering on its image.  By definition, the (projective
model of) Anti de Sitter space is the image of $X$:
\[
   AdS_3:=\pi(X)=\{[x]\in\mpi^3|\langle x,x\rangle<0\}\,.
\]
Since the covering transformation of $\pi$ preserves the metric, a
Lorentzian metric is defined on $AdS_3$. 
It is a geodesically complete Lorentzian manifold of constant curvature $-1$.
(Some authors define the AdS space as the double cover of the
$AdS_3$ space defined here, but this only introduces minor differences in the
notations.)

Notice that $AdS_3$ is an open domain in $\mpi^3$ whose boundary is
the projective quadric
\[
   \partial_\infty AdS_3:=\{x\in\mpi^3|\langle x,x\rangle=0\}~.
\]
 This quadric is a doubly ruled surface: this precisely means that there
 are two foliations $\mathcal F_l$ and $\mathcal F_r$ on
 $\partial_\infty AdS_3$ whose leaves are projective lines and such
 that the intersection of a leaf $l\in\mathcal F_l$ with a leaf
 $l'\in\mathcal F_r$ is exactly one point.

Topologically $\partial_\infty AdS_3$ is a torus and it disconnects
$\mpi^3$ in two solid tori.  It is possible to orient the leaves of
each foliation $\mathcal F_l$ and $\mathcal F_r$ in such a way that if
$c_l$ and $c_r$ denote respectively the homology classes of the
oriented leaves of the two foliations then the meridian corresponding
to $AdS_3$ is homologous to $\pm(c_l+c_r)$ and the meridian
corresponding to the complement of $AdS_3$ is homologous
$\pm(c_r-c_l)$. There are two possible way to choose such
orientations.  We fix arbitrary one of these choices.  We consider on
the boundary of $AdS_3$ the orientation such that if $e_l$ is a
positive vector tangent to the left foliation at $p$ and $e_r$ is the
positive vector tangent to the right foliation then $(e_l, e_r)$ is a
positive basis of $T_p\partial_\infty AdS_3$.  Moreover we consider on
$AdS_3$ the orientation that is compatible with the orientation of the
boundary.

The space $AdS_3$ is not simply connected. Nevertheless, isometries
act transitively on the orthonormal frames.  This implies that every
Lorentzian manifold of constant curvature $-1$ is equipped with a
$(Isom_0, AdS_3)$-structure.

Geodesics in $AdS_3$ are projective lines. There is a fairly simple
way to distinguish timelike from spacelike geodesics.  In fact
timelike geodesics correspond to projective lines entirely contained
in $AdS_3$.  They are closed simple lines of length $\pi$. Lightlike
lines correspond to projective lines that are tangent to the
boundary. Finally spacelike lines correspond to projective lines that
meet the boundary in two different points. They are open geodesics of
infinite length.

As a consequence, totally geodesic planes are obtained by intersecting
$AdS_3$ with projective planes.  Still in this case there is a
topological way to distinguish spacelike planes from timelike and
lightlike planes.  Indeed lightlike planes correspond to projective
planes tangent to $\partial_\infty AdS_3$ (that intersects the boundary
along two leaves). Timelike planes are topologically Moebius bands
(they cut $\partial_\infty AdS_3$ along a meridian of the exterior of
$AdS_3$).  Finally spacelike planes are compression disks (and they
cut $\partial_\infty AdS_3$ along a meridian of $AdS_3$).

Notice that any spacelike plane $P$ can be oriented by requiring that
its trace at infinity with the induced orientation is homologous to
$c_l+c_r$. This is called the positive orientation of $P$.  We fix the
following time-orientation  on $AdS_3$: a
timelike vector $v$ at some point $p\in AdS_3$ is future-pointing, if
it induces on the spacelike plane $P$ through $p$ orthogonal to $v$
the positive orientation.

Intrinsically a spacelike plane is isometric to $\mathbb H^2$. Indeed
it is a simply connected geodesically complete surface of constant
curvature $-1$. Moreover it can be shown that every isometry between $\mathbb
H^2$ and a plane $P_0$ extends to a projective map
\[
        r:\mpi^2\rightarrow\mpi^3  
 \]
 (where we are using the projective model of $\mathbb H^2$).
 
 In particular $r$ identifies $\partial\mh^2$ with $\partial_\infty P_0$. 
 We can consider the maps
 \[
     \partial \mh^2\rightarrow\mathcal F_l\qquad\partial\mh^2\rightarrow\mathcal F_r
 \]
 that associated to a point $p\in\partial\mh^2=\partial_\infty P_0$
 the left and the right leaves through it.  By transversality both
 maps are local homeomorphisms and for homological reasons they have
 degree one, so these maps are homeomorphisms. This precisely means
 that every leaf of $\mathcal F_l$ (resp. $\mathcal F_r$) meets $\partial\mh^2$
 exactly in one point.
 
We fix once and for all an isometric totally geodesic embedding
\[
   r_0:\mh^2\rightarrow P_0
\]
and we consider the induced identification
$\partial_\infty AdS_3$ and $\partial\mh^2\times\partial\mh^2$. Namely
any point $p\in\partial_\infty Ad_3$ is identified to the pair
$(x_l(p), x_r(p))$ where $x_l(p)$ (resp. $x_r(p)$) is the intersection
of the left (resp. right) leaf through $p$ with $\partial\mh^2=\partial_\infty P_0$.

Since isometries of $AdS_3$ are projective maps that leave
$\partial_\infty AdS_3$ invariant, then preserve the double
ruling of $\partial_\infty AdS_3$. In particular the action of
$Isom_0$ on $\partial\mh^2\times\partial\mh^2$ is diagonal: for every
$f\in Isom_0$ we have $f(x,y)=(a_l(f)(x), a_r(f)(x))$ where $a_l(f)$
and $a_r(f)$ are homeomorphism of $\partial\mh^2$.

\begin{lem}\cite{mess-notes}
The maps $a_l(f)$ and $a_r(f)$  extends to isometries of
$\mh^2$. In particular $a_l(f),a_r(f)\in PSL_2(\mr)$.
\end{lem}

By this lemma a homomorphism
\[
   a: Isom_0\ni f\mapsto (a_l(f),a_r(f))\in PSL_2(\mr)^2
\]
is pointed out. If $a_l(f)=a_r(f)=Id$, it turns out that $f$ fixes
$\partial_\infty AdS_3$. Since $f$ is a projective map, it follows
that $f=Id$. Thus $a$ is injective. Since both $Isom_0$ and $PSL_2(\mr)^2$
have dimension $6$, it follows that
$a$ is also surjective, thus it is an isomorphism.

From now on, we use the map $a$ to state an identification between
$PSL_2(\mr)\times PSL_2(\mr)$ and $Isom_0$.

\begin{remark}
The identification between $\partial_\infty AdS_3$ with
$\partial\mh^2\times\partial\mh^2$ is well-defined once we fix a
totally geodesic embedding $ r_0:\mh^2\rightarrow AdS_3$.

The map $r_0$ is unique up to post-composition with isometries of $AdS_3$.
It follows that the identifications between
$\partial_\infty AdS_3$ and $\partial\mh^2\times\partial \mh^2$ and between
$Isom_0$ and $PSL_2(\mr)\times PSL_2(\mr)$ are uniquely determined up to isometries of
$AdS_3$.
\end{remark}

Spacelike planes are determined by their
intersection with $\partial_\infty AdS_3=\partial\mh^2\times\partial\mh^2$.
By our description it turns out that the trace at infinity
of any spacelike plane is the graph of of some map $A\in PSL_2(\mr)$.

Indeed by definition the trace at infinity of our fixed plane $P_0$
corresponds to the diagonal of $\partial\mh^2\times\partial\mh^2$.  If
$P$ is any other plane, there is an isometry $f$ of $AdS_3$ such that
$f(P_0)=P$. Thus by definition 
$\partial_\infty P=\{ (a_l(f)x,a_r(f)x)|x\in\partial\mh^2\}$. By setting
$y=a_l(f)x$ we can also write
\[
\partial_\infty P=\{(y,a_r(f)a_l(f)^{-1} y)|y\in\partial\mh^2\}\,.
\]
that is, $\partial_\infty P$ is the graph of $a_r(f)a_l(f)^{-1}$.

Eventually spacelike planes are parameterized by elements in $PSL_2(\mr)$.
Given $A\in PSL_2(\mr)$ we denote by $P_A$ the plane whose
trace at infinity is the graph of $A$.

By this description, it is clear that given two planes $P,Q$
there is a unique $A$ in $PSL_2(\R)$ such that $(1,A)\cdot
P=Q$. Moreover the stabilizer of every plane is conjugated to the
diagonal subgroup into $PSL_2(\R)\times PSL_2(\R)$.

In what follows we will also use the following 
\begin{lem}
The map
\[
   \partial\mh^2\ni x\mapsto (x, Ax)\in\partial P_A
\]
extends uniquely to an isometry $r_A:\mathbb H^2\rightarrow P_A$.
\end{lem}
\begin{proof}
It is sufficient to define $r_A=(1,A)\circ r_0$.
\end{proof}

\subsubsection*{Bending angles between spacelike planes}
In Lorentzian geomety there is a natural definition of angles between
future-oriented timelike vectors. Indeed the set of future-oriented 
unit timelike tangent vectors at a point $p$ of some Lorentzian manifold $X$, say $H_p$,
is isometric to $\mathbb H^2$.
If $v,w$ lie in $H_p$ we can define the angle between $v$ and $w$ as the distance
in $H_p$ of $v$ and $w$.
By a classical formula of hyperbolic geometry, it turns out that this angle is
\[
    \mathrm{cosh}^{-1}|\langle v,w\rangle|\,.
\]

Notice that the definition is quite similar to the classical definition of angles
in Riemannian geometry, the main difference being that the angle is a well-defined number
in $[0,+\infty)$.

If $P$ and $Q$ are spacelike totally geodesic planes in $AdS_3$ 
meeting along a geodesic $l$, then their
future-oriented unit normal vector fields are parallel along $l$.
Thus we can define the bending angle between $P$ and $Q$ as the angle
between those vector fields.

If $l$ is oriented we can also define a signed bending angle between
$P,Q$.  Indeed, given a point $p\in l$, let $v,u,w\in TAdS_3$ be
respectively the positive unit tangent vector along $l$, the
future-pointing unit normal vector of $P$ and the future pointing unit
normal vector of $Q$. We say that the angle between $P$ and $Q$ is
positive if the vectors $v,u,w$ form a positive basis of $AdS_3$.  It
can be shown that the signed angle is
\[
  \alpha(P,Q)=\textrm{sinh}^{-1}\omega_p(v,w,u)
\]
where $\omega$ is the volume-form on $AdS_3$.
Notice that by definition the angle is skew-symmetric and it depends on the choice of the orientation
of $l$.

\subsubsection*{The conformal structure on $\partial_\infty AdS_3$}

Let us identify $\partial_\infty AdS_3$ with $\partial\mh^2\times\partial\mh^2$.  
If $\theta$ and $\phi$ denote positive-oriented  parameters on
each copy of $\partial\mh^2$, then we can consider the Lorentzian metric
$\eta= d\theta d\phi$ on $\partial_\infty AdS_3$.  The conformal class
of $\eta$ is independent of the choice of coordinates and  the
group $PSL_2(\R)^2$ acts conformally on $\partial\mh^2\times\partial\mh^2$.

Intrinsically the conformal structure on $\partial_\infty AdS_3$ is
characterized by the fact that isotropic directions are tangent to the
leaves of the double ruling on $\partial_\infty AdS_3$.  Indeed it can
be shown that the conformal structure on $\partial_\infty AdS_3$ is
asymptotic in the following sense: if $p_n$ is a sequence in $AdS_3$
converging to $p\in\partial_\infty AdS_3$ and $v_n\in T_{p_n}AdS_3$ is
a sequence of timelike vectors converging to $v\in
T_p(\partial_\infty AdS_3)$, then $v$ is not spacelike.

\section{Earthquakes and bent surfaces in $AdS_3$}

\subsection{Bent surfaces in $AdS_3$} 

An embedded topological surface $S\subset AdS_3$ is achronal if
geodesics joining two points of $S$ are not timelike\footnote{This
  definition is different from the standard definition of achronality in
  Lorentzian geometry -- in fact achronality does not makes sense in
  $AdS_3$ since the future of every point is the whole $AdS_3$.  On
  the other hand, if $S$ is achronal in this sense, then it is
  achronal in the standard sense in some neighbourhood}.

If $S$ is achronal then every small  neighbourhood $U$ of any
point $p\in S$ is disconnected by $S$ in two components: one is the
future of $S$ in $U$ and the other is the past of $S$ in $U$. We say
that $S$ is \emph{past convex} (resp. \emph{future convex}) if there 
is a family of neighbourhoods $\{U_i\}$ that cover $S$ and such that
for every $p,q\in S\cap U_i$ the
geodesic segment joining $p$ to $q$ does not contain points in the
future (resp. in the past)  of $S$ in $U_i$.

The surface $S$ is past convex iff for every point $p\in S$ there is a
spacelike plane $P$ such that $P\cap S$ is a convex set of $P$, and
planes obtained by moving $P$ slightly in the future do not meet
$U_i\cap S$. We say that $P$
is a support plane for $S$ in $p$. Notice that in general there are
several support planes passing through a point $p\in S$.

 \begin{remark}\label{2planes:rk}
Let $l=P\cap Q$ be oriented so that $\alpha(P,Q)>0$.  
Consider the component, of $P\setminus l$, say $P_r$, on the right side of $l$
and the component of $Q\setminus l$, say $Q_l$,   on the left side of $l$.
Then the surface
\[
   S=P_r\cup l\cup Q_l
\]   
is past convex.
In fact  we need to check the convexity only around points on $l$.
But if we slightly move $P$ in the future then the intersection with both
$P_r$ and $Q_l$ is empty.
\end{remark}

\begin{center}
\begin{figure}
\input MBH_bend.pstex_t
\end{figure}
\end{center}

A \emph{past bent surface} (resp. \emph{future bent surface}) in
$AdS_3$ is a topological embedding
\[
   b:\mathcal H\rightarrow AdS_3
\]
where $\mathcal H$ is an open convex subset of $\mh^2$ with geodesic boundary and
$b$ satisfies the following conditions:
\begin{itemize}
\item There is a geodesic lamination $L$ of $\mathcal H$ such that the
  restriction of $b$ on each connected component of $\mathcal
  H\setminus L$ is isometric and totally geodesic.
\item Each leaf of $L$ is isometrically sent  to a geodesic of $AdS_3$.
\item The image of $b$ is past convex (resp. future convex).
 \end{itemize}

\begin{remark}
A natural question is whether the map $b$ extends to the bounday.
If a boundary component $l$ of $\mathcal H$ 
is a boundary component of some stratum of $L$ then
it is clear that $b$ extends on $l$.

Instead, if there is a sequence of leaves $l_n\in L$ converging to $l$, then 
there are several possibilities:
\begin{enumerate}
\item $b(l_n)$ converges to a spacelike geodesic in $AdS_3$;
\item $b(l_n)$ converges to a point in $\partial_\infty AdS_3$;
\item $b(l_n)$ converges to a lightlike segment in $\partial_\infty AdS_3$.
\end{enumerate}
It is then clear that not in all cases the map $b$ can be extended on the boundary.
\end{remark}

There is a transverse measure on $L$ that encodes the amount of
bending along $L$.  When $L$ is locally finite, there is a fairly
simple way to describe this measure.  Given a leaf $l$ there are
exactly two regions $F,F'$ bounded by $l$, then the weight of $l$ is
simply the bending angle between the spacelike planes containing
$b(F)$ and $b(F')$.  In the general case the measure is defined by an
approximation argument using the fact that if $P,Q,R$ are spacelike
planes such that the intersection $P\cap Q$ lies \emph{above} $R$,
then the bending angle between $P$ and $Q$ is greater than the sum of
bending angles between $P$ and $R$ and between $R$ and $Q$ (see \cite{benedetti-bonsante} to
check details). 

If $c:[0,1]\rightarrow\mathcal H$ is a path transverse to $L$
then for every partition $I=(t_0=0<t_1<\ldots<t_k=1)$ one defines
$\mu(c;I)$ as the sum of the bending angles between
support planes at $b(c(t_i))$ and $b(c(t_{i+1})$.
If $I'$ is finer that $I$, the property expressed above
shows that $\mu(c,I')\leq\mu(c,I)$, so  the mass of $c$ is defined as
\[
   \mu(c)=\inf_{I}\mu(c;I)=\lim_{|I|\rightarrow 0}\mu(c;I)~. 
\]

\subsection{From earthquakes to bent surfaces}\label{earthtobend:sec}

Given two metrics $\eta_l$ and $\eta_r$ in $\mathcal T_{g,n}$ let
$h_l,h_r:\pi_1(\Sigma)\rightarrow PSL_2(\R)$ be the corresponding
holonomies.  We consider the isometric action of $\pi_1(\Sigma)$ on
$AdS_3$ given by the product holonomy
\[
   (h_l,h_r):\pi_1(\Sigma)\rightarrow PSL_2(\R)\times PSL_2(\R)~.
\]   

In this section we will associate to every right earthquake transforming $\Sigma_{\eta_l}$ into
$\Sigma_{\eta_r}$ a past bent surface 
that is invariant under the representation $(h_l,h_r)$
and we will show that this bent surface is sufficient to recover the earthquake.

Take a measured geodesic lamination  $\lambda$  on
$\Sigma_{\eta_l}$ such that the right earthquake along $\lambda$ transforms
$\Sigma_{\eta_l}$ to $\Sigma_{\eta_r}$:
\[
    E_\lambda^r:  \Sigma_{\eta_l}\rightarrow \Sigma_{\eta_r}\,.
\]
The lifting of $E_\lambda^r$ to the universal covering is a map
\[
    \tilde E:\mathcal H_l\rightarrow\mathbb H^2
 \]
 that satisfies the following properties, already mentioned in the previous section.
  \begin{itemize}
 \item The image of $\tilde E$ is the universal cover $\mathcal H_r$ of $\Sigma_{\eta_r}$.
 \item $\tilde E\circ h_l(\gamma)=h_r(\gamma)\circ \tilde E$.
 \item The image of $\tilde E$ is the convex hull of the limit set of $h_r$.
 \item For every component $F$ of $\mathbb H^2\setminus\tilde\lambda$ there is an element
 $A=A(F)$ in $PSL_2(\R)$ such that $\tilde E|_F=A|_F$.
 \item $A(h_l(\gamma)(F))=h_r(\gamma)\circ A(F)$.
 \item If $F,F'$ are two components of $\mathcal
   H\setminus\tilde\lambda$ then the comparison isometry $B^*=A(F)^{-1}\circ
   A(F')$ is a hyperbolic transformation whose axis separates $F$
   from $F'$.  If the axis $l$ of $B^*$ is oriented from the repulsive
   fixed point towards the attractive fixed point, then $F'$ is on the
   left side of $l$ whereas $F$ is on the right side.
 \item If $F$ and $F'$ are adjacent then the axis of $B^*$ is the common
   edge $e$ and the translation length of $B^*$ is the weight of $e$.
 \end{itemize}

Given any component $F$ in $\mathcal H_l\setminus\tilde\lambda$ let us take the set
of its ideal vertices $\{x_i\}\subset\partial\mh^2$.
Setting $A=A(F)$, we can consider on the plane $P_A$ the convex hull of the
set $\{(x_i, Ax_i)\}$, that is a spacelike geodesic polygon in $AdS_3$, say  $K(F)$.

Let $S$ be the closure of the union of all $K(F)$'s.

\begin{prop}\label{earthtobend:prop}
$S$ is a future convex bent surface in $AdS_3$ that is invariant under
  the action of $\pi_1(S)$.

Moreover if $\mathcal H$ denotes the 
universal covering of  $E^r_{\lambda/2}(\Sigma_{\eta_l})$
then there is a bending map
\[
  \iota:\mathcal H\rightarrow S 
\]
that is equivariant under the $\pi_1(\Sigma)$-action.

The bending lamination associated to $\iota$ is the image throught the
earthquake map $E:\mathcal H_l\rightarrow\mathcal H$ of
the lamination $\tilde\lambda/2$.
\end{prop}

To prove Proposition \ref{earthtobend:prop} we need the following elementary facts of AdS geometry,
the proofs can be found in \cite{benedetti-bonsante}.

\begin{lem}\label{comp:lm}
 Let $l$ be a complete geodesic line in $AdS_3$ with endpoints 
 $p=(x,y)$ and $q=(x',y')$. Let $s_l, s_r$ be respectively the geodesics
 of $\mathbb H^2$ with end-points $x,x'$ and $y,y'$.
  
The connected component of the stabilizer of $l$ in $PSL_2(\R)\times PSL_2(\R)$ is the set
of pairs $(A,B)$ where $A$ is a hyperbolic transformation with axis $s_l$ and $B$ is
a hyperbolic transformation with axis $s_r$.

Let us orient $l$, $s_l$ and $s_r$ in such a way that the
corresponding starting points are respectively $(x,y), x$ and $y$.
Given a transformation with axis $A$ (resp. $B$) let $t(A)\in\mathbb R$
(resp. $t(B)$) denote the signed translation length ($t(A)$ is positive if
$x$ is the repulsive fixed point, negative otherwise).
Then the transformation $(A,B)$ acts as a translation on $l$ of factor $(t(A)+t(B))/2$.
The rotation angle of $(A,B)$ along $l$ is $(t(B)-t(A))/2$
\end{lem}

\begin{defi}
If $(A,B)$ preserves $l$, the rotation angle of $(A,B)$ along $l$ is
the signed  bending angle formed by a spacelike 
plane $P$ containing $l$ with its image $(A,B)\cdot P$.
\end{defi}

We prove now Proposition \ref{earthtobend:prop}.

\begin{proof}[Proof of Proposition \ref{earthtobend:prop}]
We prove the statement assuming that $\lambda$ is locally finite. The general case will
follow by an approximation argument.

For each face of $\mathcal H_l\setminus\tilde\lambda$ let $r_F:\mh^2\rightarrow AdS_3$ 
be the isometric embedding whose trace at infinity is the map
\[
     x\mapsto (x, A(F)x)
\]
Clearly we have that $K(F)=r_F(F)$.

Give two strata $F,F'$ we have that $r_{F'}=(1,B)\circ r_F$ where 
$B=A(F')\circ A(F)^{-1}$.
Thus $K(F')$ is obtained by applying the transformation $(1,B)$ to 
$r_F(F')$.

 Notice that $B=A(F)B^*A(F)^{-1}$ where
$B^*=A(F)^{-1}A(F')$ is the comparison isometry.

Let $l$ be the image through $r_F$ of the axis of $B^*$,
that is the geodesic with end-points $p_-=(x_-(B^*), A(F)x_-(B^*))=(x_-(B^*), x_-(B))$
and $p_+=(x_+(B^*), A(F)x_+(B^*))=(x_+(B^*),x_+(B))$.
Notice that both $p_-$ and $p_+$ are fixed by $(1,B)$.
Thus $l$ is left invariant by $(1,B)$.

Now if $P_l$ and $P_r$ denote the half-planes bounded by $l$
on $P_{A(F)}$, then $K(F)=r_F(F)$ is contained in 
$P_r$ whereas $r_F(F')$ is contained in $P_l$.

We can conclude that
\begin{equation}\label{cont:eq}
K(F)\cup K(F')\subset P_r\cup (1,B)P_l\ .
\end{equation}
This shows that if $F$ and $F'$ are not adjacent, then
$K(F)$ and $K(F')$ are disjoint. When $F$ and $F'$ meet
along a line, this line is the axis of $B$, so $K(F)$ and $K(F')$
meet along the line $l$. In this case
by Lemma \ref{comp:lm}, the bending angles formed along $l$
between $K(F)$ and $K(F')$ is equal to $t(B)/2=t(B^*)/2$, that is to one half
the mass of the line $F\cap F'$ 

Since $P_r\cup (1,B) P_l$ is achronal, it turns out that 
$S$ is achronal.

Notice that  we have that
\[
   K(h_l(\gamma)(F))=(h_l(\gamma), h_r(\gamma))K(F)
\]   
thus $S$ is invariant.

In order to show that $S$ is a past bent surface, we need
to construct the bending map.

We could try to glue the maps $r_F$.
However if $F$ and $F'$ are adjacent, for $p$ in $F\cap F'$ we have
\[ 
    r_{F'}(p)=(1,B)r_F(p)~.
\]
Notice that both $r_F(p)$ and $r_{F'}(p)$ are contained in the geodesic
$l$ described above. On the other hand, by Lemma \ref{comp:lm}
the transformation  $(1,B)$ acts by a translation of factor
$t(B)/2$ on $l$.
Thus the maps $r_F$ do not glue to an isometric map from $\mathcal H_l$
into $AdS_3$.  

On the other hand, these maps can be glued if each
 component of $\mathcal H_l\setminus\tilde\lambda$ is identified to
 the adjacent components through a right translation of length equal to
 the mass of the corresponding edge divided by $2$.  This shows that
 the maps $r_F$ induce a continuous isometric identification
  \[
  \mathcal H\rightarrow S~.
  \]
  
 To conclude we have to prove that $S$ is past convex.  It is
 sufficient to show that $S$ is convex at each point in $K(F)\cap
 K(F')$ where $F$ and $F'$ are two adjacent components of $\mathcal
 H_l\setminus\lambda$. On the other hand by (\ref{cont:eq}) we have that
 $K(F)\cup K(F')$ is contained in $P_r\cup (1,B)P_l$. By Lemma
 \ref{comp:lm} the angle formed between $P_{A(F)}$ and $(1,B)P_{A(F)}$
 is positive (with respect to the natural orientation of $l$).  Thus
 Remark \ref{2planes:rk} shows that $P_r\cup (1,B)P_l$ is past convex
 and this concludes the proof.
 \end{proof}

\begin{remark}
The surface $S$ is well-defined up to post-composition with an
isometry of $AdS_3$.  Notice that different earthquakes produce
different bent surfaces. This depend on the fact that the shearing
lamination of the earthquake
$E^r:\Sigma_{\eta_l}\rightarrow\Sigma_{\eta_r}$ is explicitly related
to the bending lamination of $S$.
\end{remark}

It could be proved that every bent convex surface 
can be associated to an earthquake. More precisely the following statement holds.

\begin{prop} \label{pr:57}
Consider an equivariant bent surface
\[
    b:\mathcal H\rightarrow AdS_3
\]
such that $S_+=\mathcal H/(h_l,h_r)$ is an admissible surface.
Let $\lambda_+$ be the lamination on $S$ corresponding to the bending
lamination on $H$.  Then
\[
 E^r_\lambda(S_+)=S_l~, \qquad E^l_\lambda(S_+)=S_r~.
\]
\end{prop}

Since this proposition is not strictly necessary for the proof of Theorem \ref{tm:earthquake}
we skip the proof referring to \cite{benedetti-bonsante} the interested reader.

\section{Bent surfaces in $AdS_3$ and achronal meridians in $\partial_\infty AdS_3$}

This section analyses the relationship between bent surfaces in the AdS space and
curves of a certain type -- called achronal meridians -- arising as their boundary
at infinity. 

\subsection{Achronal meridians as graph of cyclic-order  preserving maps of $\partial\mh^2$} \label{achronal:sec}

In this section it will be convenient to fix a universal covering of
$\partial\mathbb H^2\times\partial\mathbb H^2$.  In particular we fix
a point $q_0$ in $\mathbb H^2$ and we consider the visual angle on
$\partial\mh^2$ with respect to $q_0$.  This gives a natural covering
map $p:\mathbb R\rightarrow\partial\mh^2$.  Clearly we can consider
the product covering $\mathbb
R^2\rightarrow\partial\mh^2\times\partial\mh^2$ sending
$(\theta,\phi)$ to $(p(\theta), p(\phi))$.
For notation convenience we  slightly modify this covering
by considering the map
\[
  \mathbb R^2\ni (x,y)\mapsto (p(2\pi x), p(2\pi
  y))\in\partial\mh^2\times\partial\mh^2\,.
  \] 
In this way the covering transformations are translations with integer
coordinates.

By definition spacelike curves on $\partial_\infty AdS_3$ correspond
to curves $(x(t), y(t))$ such that
$x'(t)y'(t)>0$. They are locally graphs of orientation
preserving maps between two open intervals of $\partial\mh^2$.  Thus, the
lifting on the universal covering $\R^2$ of a spacelike curve $c$ is
the graph of a strictly increasing function
\[
    f:\mathbb R\rightarrow\mathbb R
\]
such that $f(x+n)=f(x)+m$ for some $n,m\in\Z$ depending on the
homology class of $c$.

Spacelike meridians in $\partial_\infty AdS_3$ are graphs of
orientation-preserving diffeomorphisms of $\partial\mh^2$, since their
liftings correspond to graphs of orientation preserving
diffeomorphisms $f:\mathbb R\rightarrow\R$ such that $f(x+1)=f(x)+1$.

Limit of spacelike meridians are locally achronal meridian. A meridian
is {\it locally achronal} if for every $p\in c$ there is a neighborhood
$U\subset\partial_\infty AdS_3$ such that no pair of points $q,r\in
c\cap U$ are related by a timelike arc in $U$.  It can be shown that
locally achronal meridians correspond to monotonically increasing
(possibly discontinuous) functions
\[
   f:\R\rightarrow \R
\]
such that $f(x+1)=f(x)+1$.  Indeed given such a function, we can
consider the subset in $\R^2$
\[
  G_f=\{ (x,y)| \lim_{t\rightarrow x_-} f(t)\leq
  y\leq\lim_{t\rightarrow x_+} f(t)\}\,.
\]
\begin{lem}
$G_f$ is a connected embedded curve in $\R^2$
\end{lem}
\begin{proof}
If $f_n$ is a sequence of continuous monotonically increasing
functions approximating $f$ point-wise,  the
length of the graph of $f_n$ on some interval $[a,b]$ is bounded by
$(b-a)+(f_n(a)-f_n(b))$ so it is uniformly bounded. Since
$graph(f_n|_{[a,b]})$ stays in some compact set of $\mathbb R^2$ it
converges to a topological curve. Such a curve coincides with
$G_f|_{[a,b]}$.

The fact that $G_f$ is embedded is due to the fact that every point of
$G_f$ disconnects $G_f$ in exactly two components.
\end{proof}

Given an increasing function $f$, the set $G_f$ projects to an embedded
closed curve in $\partial\mh^2$ provided that $f$ is not constant on
some interval of length bigger than $1$ or discontinuity points with
jumps bigger than $1$.  On the other hand by our assumption
$f(x+1)=f(x)+1$, it is easy to see that this can happen if and only
if up to some translation we have $f(t)=[t+c]+c'$ for some constants
$c,c'$.  In both these cases $G_f$ projects in the union of two leaves
in $\partial\mh^2$.

In all the other cases $G_f$ projects to a locally achronal meridian
in $\partial\mh^2\times\partial\mh^2$. Conversely every locally
achronal curve arises in this way. Namely, given an achronal meridian
$C$ we define $f:\mathbb R\rightarrow\mathbb R$ by setting
$f(x)=\sup\{y|(x,y)\in\tilde C\}$ where $\tilde C$ is a component of
the pre-image of $C$ in $\mathbb R^2$.  Since $C$ meets every leaf,
the map is well-defined.  Since $C$ is a meridian,
$f(x+1)=f(x)+1$.  Finally, since $C$ is achronal $f$ turns out to be
increasing and $C$ coincides with the projection of $G_f$.

With some abuse we call $G_f$ the graph of the
function $f$ (notice that $G_f$ coincides with the standard graph when
$f$ is continuous).

We collect some facts about locally achronal meridians that will be useful in what follows.

\begin{lem}
If $C$ is an achronal meridian then 
for every point $p\in\partial_\infty AdS_3\setminus C$ 
there is a projective plane $P$ containing $p$
whose intersection with $AdS_3$ is spacelike and
such that $C\cap P=\emptyset$.
\end{lem}

\begin{proof}
Let $f:\mathbb R\rightarrow\mathbb R$ be the increasing function such that $C$
is the projection of $G_f$.
Acting by isometries  on $C$ we can 
suppose that  $f$ is continuous at $0$ and $1/2$ 
and  $f(0)=0$ (so $f(1)=1$) and $f(1/2)=1/2$.
Notice that $G_f\cap [0,1]^2$ is contained in the two squares 
$Q_1=[0,1/2]\times [0,1/2]$ and $Q_2=[1/2,1]\times[1/2,1]$.
Moreover, by our assumption on the continuity
points $(1/2,0)$, $(0,1/2)$ , $(1,1/2)$ and $(1/2,1)$do not  lie in $G_f$.

Thus the line of equation $y=x+1/2$ is disjoint from $G_f$.
The projection $l$ of such a line on
$\partial\mh^2\times\partial\mh^2$ is the graph of the trace at
infinity of the rotation of angle $\pi$ about the
point $q_0$.  Thus there is a spacelike plane $P$ such that $l=P\cap
\partial_\infty AdS_3$.  It follows that $P\cap C=\emptyset$.
 \end{proof}

Since $C$ does not intersect $P$ we can consider the convex hull $K$
of $C$ in the affine chart $\mathbb R^3=\mpi^3\setminus P$. It is easy
to see that $K$ does not depend on the plane $P$ (this because the
change of chart map between $\mathbb R^3\setminus P$ and $\mathbb
R^3\setminus Q$ sends compact convex sets disjoint from $Q$ 
into compact convex sets).

\begin{lem}
$K$  is contained in $\overline{AdS_3}$.
More precisely 
\begin{enumerate}
\item the interior of $K$ is contained in $AdS_3$
\item the intersection of the boundary of $K$ with $\partial_\infty
  AdS_3$ is $C$.
\end{enumerate} 
\end{lem}

\begin{proof}
Given a point $p=(x,y)$ in $\partial AdS_3\setminus C$, we claim that there
exists a spacelike plane passing through $p$ that does not intersect
$C$. As a consequence we have that $K\cap\partial_\infty AdS_3=C$.  In
particular $K$ is contained in the closure of one component of
$\mathbb R^3\setminus \partial_\infty AdS_3$. Since the curve $C$ is
not trivial in $\mathbb R^3\setminus AdS_3$, $K$ must be contained in
$AdS_3$.

We now prove the claim.  Let us consider a timelike plane, say $Q$,
through $p$.  For homological reasons $Q$ must intersect $C$ in two
points $q=(u,v),q'=(u',v')$.  

Since $Q\cap\partial_\infty AdS_3$ is the graph of an orientation
reversing diffeomorphism of $\partial\mh^2$, we can fix $q,q'$ so that
$x$ is contained in the positive segment $[u,u']$ whereas $y$ is in
the positive segment $[v',v]$.

Up to the action of $PSL_2(\mr)\times PSL_2(\mr)$ we can also suppose that
the points 
\[
  \hat q=(0,0)\qquad \hat p=(1/3,2/3)\qquad \hat q'=(2/3,1/3)
\]
project respectively to $q,p,q'$.  Let $\tilde C$ be the lifting of $C$
passing through $(0,0)$.  We have that $\tilde C\cap[0,1]^2$ is
contained in $R=[0,2/3]\times[0,1/3]\cup[2/3,1]\times[1/3,1]$.
Consider the line $\hat l$ of equation $y=x+1/3$.  There is a
spacelike plane $P$ such that $P\cap\partial_\infty AdS_3$ is the
projection of $\hat l$.  If the points $(0,1/3), (2/3,0), (2/3,1)$ do not
lie on $\tilde C$, the plane $P$ has the required
property. Otherwise, since $\hat l$ does not disconnect $\tilde C$,
the plane $P$ is a support plane for $K$.  Notice that in this case
$K\cap P$ is the convex hull of $C\cap P$. Since $\partial_\infty
AdS_3\cap P$ is strictly convex, we have that $p$ does not lie on
$K\cap P$, so $p\notin K$.
\end{proof}

\begin{remark}
Support planes of $K$ cannot be timelike, indeed for homological
reasons the transverse intersection of $C$ with a timelike plane is
not empty.

$K$ is a plane if and only if $C$ lies in some projective plane, otherwise it is
a closed ball (in $\mpi^3$).  The boundary of $K$ in $AdS_3$ has two
connected components. By the remark above both components are achronal
surfaces. More precisely one component is past convex 
and the other is future convex.
\end{remark}

The \emph{upper boundary} of $K$ -- denoted by $\partial_+K$ -- is the
past convex component of $\partial K$. Analogously the \emph{lower
  boundary} of $K$ -- denoted by $\partial_-K$ -- is the future convex
component of $\partial K$.

We say that a support plane of $K$ is an upper (resp. lower) 
support plane if it is a support for the upper boundary
(resp. lower boundary).

\begin{remark}\label{support:rk}
Let $f:\R\rightarrow\R$ be an increasing function
such that  $C$ is the projection of $G_f$.
Given $A\in PSL_2(\R)$, the following conditions are equivalent:
\begin{enumerate}
\item $P_A$ is a upper (resp. lower) support plane for $C$.
\item There is a lifting of $A$ to $\mathbb R$ such that
$f_A=\tilde A^{-1}\circ f$ satisfies $f_A(x)\leq x$ (resp. $f_A(x)\geq x$) 
and admits two fixed points on $[0,1)$.
\end{enumerate}
\end{remark}

If $P$ is a support plane for $K$, the intersection $P\cap K$ is the convex hull
of $P\cap C$.

If $P$ is a spacelike support plane the intersection of $P$ with $\partial K$ is
either a geodesic line or a  hyperbolic ideal polygon.

If $P$ is a lightlike support plane, it is tangent to $\partial_\infty
AdS_3$ at some point $p$ in $C$.  Moreover since $P$ meets the
boundary of $K$, there are points $q,r$ of $C$ lying respectively on
the left and the right leaf through $p$. The intersection of $P$ with
$C$ is a lightlike triangle with two ideal edges (that means that two
edges are segments of leaves of the double foliation of
$\partial_\infty AdS_3$).

It turns out that each boundary component of $K$ is the union of a
\emph{spacelike} region formed by the set of points admitting only
spacelike supports planes and some ideal lightlike triangles.

The spacelike part is a union of space-like geodesics and of ideal hyperbolic
polygons. 
The boundary of the spacelike part in $\partial_+K$ is the union of the
spacelike edges of the ideal lightlike triangles contained in $\partial_+K$.

In what follows we will need the following technical fact.

\begin{lem}\label{intersection:lm}
Let $C$ be an achronal meridian and let $K$ denote its convex hull.
If $P$ and $Q$ are spacelike upper support planes, then they intersect along a line
$l$. Moreover if $l$ is oriented so that $\alpha(P,Q)>0$ then $P\cap K$ is contained in
the right side of $l$ in $P$ and $Q\cap K$ is contained in the left side  of $l$ in $Q$.
\end{lem}
\begin{proof}
Suppose that two upper support planes $P,Q$ are disjoint in
$AdS_3$. Then there are planes $P',Q'$ obtained by slightly moving in
the future of $P$ and $Q$ respectively such that $P'\cap K=Q'\cap
K=\emptyset$ and we can moreover suppose that $P'\cap Q'= P'\cap
Q=Q'\cap P=\emptyset$.  Notice that $AdS_3\setminus (P'\cup Q')$ is
the disjoint union of two cylinders and $P$ and $Q$ lie in different
components. Now $K$ is a connected set in $AdS_3\setminus (P'\cup Q')$
that contains a point on $P$ and a point on $Q$, so it intersects the
two components of $AdS_3\setminus (P'\cup Q')$ and this gives a
contradiction.

To conclude the proof of the Lemma, it is sufficient to notice that
given $p\in P$ and $q\in Q$, the only possibility that the geodesic
segment $[p,q]$ intersects neither $P'$ nor $Q'$ is that $p\in P_r$
and $q\in Q_l$.
\end{proof}

\begin{remark}\label{intersection:rk}
Let $P$ $Q$ and $l$ be as in Lemma \ref{intersection:lm}.  Let
$p_-=(x_-,y_-)$ and $p_+=(x_+,y_+)$ be respectively the starting and
the ending point of $l$. Let $s$ be the geodesic in $\mh^2$ with
end-points $x_-$ and $x_+$ (oriented from $x_-$ to $x_+$).  Then if
$P\cap C=\{(x_i,y_i)\}_{i\in I}$ and $Q\cap C=\{x'_j, y'_j)\}_{j\in
  J}$ points $x_i$ lie on the right side of $s$ whereas points $x'_j$
lie on the left side of $s$.  In particular $s$ weakly separates the
convex hull of points $x_i$ from the convex hull of points $x'_j$.
\end{remark}

\subsection{From bent surfaces to achronal meridians}
Let
\[
  E^r:\Sigma_{\eta_l}\rightarrow\Sigma_{\eta_r}
\]
be a right earthquake between two admissible surfaces.
In this section an achronal curve $C$ will be associated to $E^r$
and we will show that such a curve determines the earthquake.

More precisely we  consider the equivariant bent surface
\[
  b:\mathcal H\rightarrow AdS_3
\]
associated to $E^r$ in Section \ref{earthtobend:sec}.
We will construct  an achronal meridian $C$ such that 
$b(\mathcal H)$ turns out to be the spacelike region of
the convex hull of $C$.

First we consider the trace at infinity of the earthquake.
Contrary to the ``classical'' case (where the source and the target spaces are the whole $\mh^2$)
it is not true that  the map $E^r$ extends by continuity on the closure of $\mathcal H$ at infinity.
Nevertheless the map $E^r$ extends on the set of ideal points of every stratum of the fault lamination.
So we consider the following set in $\partial_\infty AdS_3$
\[
  \partial_\infty\mathcal H=\{(x,E^r(x))| x\in F\cap\partial\mathbb H^2, F\textrm{ is a stratum of the fault lamination}\}\ .
\] 

This notation is due to the fact that $\partial_\infty\mathcal H$ can be regarded as
the set of ideal points of the bent surface $b(\mathcal H)$.

More precisely, for every stratum of the fault lamination $F$,
we consider the set $\partial_\infty K(F)$ 
of ideal points of $K(F)$ where $K(F)$ is the face of $b(\mathcal H)$ corresponding to $F$.
It turns out that
\[
  \partial_\infty\mathcal H=\bigcup\partial_\infty K(F)\ .
\]

\begin{remark}\label{invariance:rk}
If $h_l$ and $h_r$ are the holonomies of $\eta_l$ and $\eta_r$ then
$\partial_\infty\mathcal H$ is invariant for the representation $(h_l,h_r)$.
\end{remark}

Given three  points $x,y,z\in\partial\mh^2$ such that $x\neq z$
we write $x\leq y\leq z$
if $y$ lies in the positive closed segment in $\partial\mh^2$ with first end-point
$x$ and second end-point $z$. We write $x<y<z$ if $x\leq y\leq z$ and $y\neq x$ and $y\neq z$.
The set $\partial_\infty \mathcal H$ satisfies the following property:

\begin{lem}\label{acronal:lm}
Given three points $p_1,p_2,p_3\in\partial_\infty\mathcal H$ such that
$p_i=(x_i,y_i)\in\partial_\infty \mathcal H$ and $x_1<x_2<x_3$
then 
$y_1< y_2< y_3$.
\end{lem}

\begin{center}
\begin{figure}
\input MBH_achron.pstex_t
\end{figure}
\end{center}

\begin{proof}
%


Let $F_i$ be the stratum of the fault lamination such that $p_i\in K(F_i)$.
We prove the statement assuming $F_1\neq F_2\neq F_3$ . The other cases are simpler and quite similar.

Up to applying some cyclic permutation of indices, we may suppose that $F_2$ separates $F_1$ from $F_3$. 
This precisely means that there are points in $\partial_\infty F_2$, say $x_1', x_1'', x_3', x_3''$ such that
\[
x_1'\leq x_1\leq x_1''~,\quad x_3'\leq x_3\leq x_3''
\]
and the positive intervals $(x_1',x_1'')$ and $(x_3', x_3'')$ are disjoint 
and do not contain points in $\partial_\infty F_2$.
Moreover by the hypothesis, either $x_1''=x_2=x_3'$ or $x_1''\leq x_2\leq x_3'$.

Now, let $A\in PSL_2(\mr)$ such that $E^r|_{F_2}=A$. By properties of the earthquake 
$A^{-1}E^r(x_i)$ is contained in $[x_i',x_i'']$ for $i=1,3$, therefore we have
\[
     A^{-1}y_1<x_2<A^{-1}y_3~.
\]
Since $A$ is an orientation preserving diffeomorphism of $\partial\mh^2$ we conclude  that
 \[
   y_1<y_2<y_3.
 \]
\end{proof}

The property expressed in Lemma \ref{acronal:lm} is shared in some weaker
form by all subsets of any achronal meridian.

\begin{lem}\label{mer:lm}
If $C$ is an achronal meridian in $\partial_\infty AdS_3$ then given three points $p_1,p_2,p_3\in C$
such that $p_i=(x_i,y_i)$ with $x_1<x_2<x_3$ then either $y_1=y_2=y_3$ or 
$y_1\leq y_2\leq y_3$.
\end{lem}

\begin{proof}
It is sufficent to notice that $C$ is the projection of some curve $G_f$ where 
$f:\mr\rightarrow\mr$ is an increasing function such that $f(x+1)=f(x)+1$.
\end{proof}

\begin{defi}
A subset of $\partial_\infty AdS_3$ that 
is not contained in the union of any left and right leaves and
satisfies the conclusion of Lemma \ref{mer:lm}
will be said to be connectible by an achronal meridian. 
\end{defi}

\begin{lem}\label{con:lm}
If $\Delta\subset \partial_\infty AdS_3$ 
is connectible by an achronal meridian
then there is an achronal meridian passing through every point of $\Delta$.

Indeed there are two extremal possible choices $C_-(\Delta), C_+(\Delta)$ such that
every other choice lies in between them.
\end{lem}

\begin{proof}
We fix the angular coordinates on $\partial\mh^2\times\partial\mh^2$ such that 
 $(0,0)$ corresponds to some point in $\Delta$.

Let $\tilde\Delta$ be  the pre-image of $\Delta$ on the open square $(0,1)\times (0,1)$
through the covering map $\mathbb R^2\rightarrow\partial\mh^2\times\partial\mh^2$
described in part \ref{achronal:sec}. 
This set has the property that if $(x,y), (x',y')\in\tilde\Delta\cap (0,1)^2$ 
and $x<x'$ then $y\leq y'$.

Thus we can define $f_-:[0,1]\rightarrow [0,1]$ by setting
\[
  f_-(0)=0~, \qquad f_-(1)=1~,\qquad
  f_-(t)=\sup\{y~|~ \exists (x,y)\in\tilde\Delta\textrm{ s.t. }  x\leq t\} 
\]
where we use the convention $\sup\emptyset=0$.
This function is clearly increasing, and we can extend
$f_-$ to an increasing function on $\mathbb R$ such that $f_-(t+1)=f_-(t)+1$

If $(x,y)\in\tilde\Delta$ then $f(x)\geq y$; on the other hand by
the property of $\tilde\Delta$ we have that
\[
\lim_{t\rightarrow x_-} f_-(t)\leq y
\]
thus $(x,y)$ is contained in $G_{f_-}$.  Finally notice that if
$(t,0)$ projects to a point of $\Delta$, then there cannot be any
point $(x,y)\in\tilde\Delta$ such that $x<t$. It follows that $f=0$ on
the interval $[0,t)$, and $(0,t)\times\{0\}\in G_f$.

$G_f$ projects to some curve $C$ in
  $\partial\mh^2\times\partial\mh^2$.  Since $C$ contains $\Delta$, it
  cannot be the union of any left and right leaves.  Thus $C$ is an
  achronal meridian containing $\Delta$.

We can also define $f_+:[0,1]\rightarrow [0,1]$ by putting
\[
   f_+(0)=0~,\qquad f_+(1)=1~,\qquad f_+(t)=\inf\{y~|~\exists
   (x,y)\in\tilde\Delta\textrm{ s.t. }  t\leq x\}~,
\]
where we use the convention that  $\inf\emptyset=1$.
The same argument used above shows that this function
is increasing and that the corresponding achronal meridian contains
$\Delta$.

Let $C$ be an achronal meridian containing
$\Delta$. There is a monotonic function $f:\mathbb R\rightarrow\mathbb R$
such that  $f(x+1)=f(x)+1$ and $C$ is the projection  of $G_f$. Clearly
we can normalize $f$ so that $f(0)=0$.

It is easy to see that on the interval $[0, 1]$ we have $f_-\leq f\leq
f_+$, so the same inequalities hold on the whole real line.
 \end{proof}

\begin{remark}\label{closure:rk}
The property to be connectible by an achronal meridian is closed. That
is, if $\Delta$ is connectible, so is $\overline{\Delta}$.
\end{remark}

\begin{remark}\label{mer:rm}
If $\Delta$ is connectible by an achronal meridian and it is connected then 
$\Delta$ is itself an achronal meridian and $C_-=C_+$.

Otherwise we can consider the region $Q(\Delta)$ obtained by
projecting to the quotient the region
\[
\tilde Q=\{(x,y)\in [0,1]^2~|~ f_-(x)\leq y\leq f_+(x)\}
\] 
It is the union of the closure of $\Delta$ and some rectangles whose
edges (that are light-like for the conformal structure of $\partial_\infty
AdS_3$) are contained in $G_{f_-}$ and $G_{f_+}$ (we consider also the
degenerate case when two opposite edges collapse to points). Two
opposite vertices of such rectangles lie in $\Delta$.

Clearly every achronal meridian containing $\Delta$
is contained in $Q$. Conversely if for every rectangle,  an achronal
arc connecting the two vertices in $\Delta$ is choosen, then the closure
of the union of such arcs is an achronal meridian.
\end{remark}

\begin{center}
\begin{figure}
\input MBH_bound.pstex_t
\end{figure}
\end{center}

\begin{remark}\label{rect:rk}
Let $\Delta=\{(x_i,y_i)\}$ be connectible by an achronal meridian.
Let $I$ be a component of $\partial\mh^2\setminus\overline{\{x_i\}}$.
Then there is a rectangle $R$ in $Q(\Delta)\setminus\Delta$ of the
form $I\times J$. Moreover this a correspondence is $1$-to-$1$.
 
 In particular given two points $(x,y), (x',y')\in\overline{\Delta}$
 with $x\neq x'$, if there is no point $p=(x'',y'')$ in $\Delta$ such
 that $x<x''<x'$ then the rectangle $[x,x']\times[y,y']$ is contained
 in $Q(\Delta)$ (here $[y,y']$ is the positive segment joining $y$ to
 $y'$ if they are different, otherwise $[y,y']=\{y\}$).
 \end{remark}

\begin{remark}\label{lightlike:rk}
Let $\Delta$ be connectible by an achronal meridian and $C_-$ be the
lower meridian through $\Delta$.  Consider a non-degenerate rectangle
in $Q(\Delta)$, say $R=[x,x']\times[y,y']$.  Recall that the lightlike
plane $P$ tangent to $\partial_\infty AdS_3$ at a point $(x',y)$ meets
$\partial_\infty AdS_3$ along the leaves through $(x',y)$.  It follows
that $P$ does not separate $C_-$, so it is a support plane for the
convex hull $K$ of $C_-$.  In particular $K\cap P$ is the lightlike
ideal triangle with vertices $(x,y), (x',y), (x',y')$.
\end{remark}

Let us come back to our original problem.  By Lemmas \ref{acronal:lm}
and \ref{con:lm} we conclude that $\partial_\infty \mathcal H$ is connectible
by an achronal meridian.  Let $C_-$ be the extremal lower meridian
passing through $\partial_\infty \mathcal H$.
\begin{remark}
Since $\partial_\infty\mathcal H$ is invariant for the representation
$(h_l,h_r)$, it is easy to check that also $C_-$ is invariant.
\end{remark}

We are going to prove that the curve $C_-$ determines the bent surface
$b(\mathcal H)$ and thus determines the earthquake $E^r$. Recall that 
$b(\mathcal H)$ is defined in part \ref{earthtobend:sec}.

\begin{prop} \label{pr:618}
The bent surface $b(\mathcal H)$ is
the spacelike part of the future boundary of the convex hull $K$ of $C_-$.
\end{prop}
\begin{proof}
For every stratum $F$ of the fault lamination the set 
 $K(F)$ is the convex hull of its ideal points, that lie in $\partial_\infty\mathcal H$.
 It follows that $b(\mathcal H)$ is contained in $K$.
 
 We claim that $b(\mathcal H)$ is contained in the upper boundary of $K$.
 Given a stratum $F$ of the fault locus, let $P_F$ denote the spacelike plane
 in $AdS_3$ containing $K(F)$. We will prove that $P_F$ is an upper support plane
 for $K$.
 
 Indeed, up to post-composition with an element in $PSL_2(\mr)$ we may
 suppose that $E^r|_F=Id$. It follows that $P_F$ is the plane $P_0$
 whose trace at infinity is the diagonal of
 $\partial\mh^2\times\partial\mh^2$.
 
 Let $F'$ be another stratum. There are two ideal points $x,x'\in
 \partial_\infty F$ such that the geodesic in $\mh^2$ with endpoints
 $x$ and $x'$ is a component of the frontier of $F$ and
 $\partial_\infty F'$ is contained in the positive interval $(x,x')$.
 The fact that $F'$ is moved on the right as seen from $F$ means that
 if $y\in\partial_\infty F'$ then either $y=E^r(y)\in\{x,x'\}$ or
 $x<E^r(y)<y$.

This shows that the pre-image of $\partial_\infty \mathcal H$ on the
square $(0,1)\times(0,1)$ (where we are assuming that the point
$(0,0)$ corresponds to an ideal point of $F$) is contained in the
triangle $\{(u,v)|v\leq u\}$.

It easily follows that $f_-(u)\leq u$ for every $u\in [0,1]$. By
Remark \ref{support:rk}, the plane $P_F$ turns out to be an upper support
plane for $K$.  This shows that $b(\mathcal H)$ is contained in the
spacelike part of $\partial_+K$.

Let $p\in\partial_+K\setminus b(\mathcal H)$. We will show that a
lightlike support plane passes through it. Suppose that $P$ is a spacelike support
plane through $p$.  Let us consider $P\cap C_-=\{(x_i,y_i)\}$, and let
$F_P$ be the convex hull in $\mh^2$ of the points $x_i$.  By Lemma
\ref{intersection:lm} and Remark \ref{intersection:rk}, we have that
$F_P$ is weakly separated by all the strata of the fault
lamination. In particular it is separated from the universal cover, $\mathcal H_l$,
of $\Sigma_\eta$.

It follows that points $x_i$ are all contained in some component $I$
of $\partial\mh^2\setminus\partial\mathcal H_l$. Thus $C\cap P_-$ is
contained in the rectangle $R$ of $Q(\partial_\infty\mathcal H)$
corresponding to the interval $I$.  $R\cap C_-$ is the union of the
two lower edges connecting vertices
$q,q'\in\overline{\partial_\infty\mathcal H}$.

Using the fact that $P$ is a spacelike support plane, it follows easily that
$P\cap C_-=\{q,q'\}$, thus $p$ lies on the geodesic with end-points $q,q'$.

By Remark \ref{lightlike:rk},
the lightlike plane tangent to $\partial_\infty AdS_3$ at the lower vertex of $R$ 
is a support plane for $K$ containing $q,q'$ and thus $p$.
It follows that $p$ is not contained in the spacelike region of $\partial_+K$.
\end{proof}

\section{The action of $\pi_1(\Sigma)$ on $\partial \mh^2\times\partial\mh^2$}

Given two elements $\eta_l,\eta_r\in\Tt_{g,n}$ with corresponding holonomies 
$h_l,h_r:\pi_1(\Sigma)\rightarrow PSL(2,\mr)$, let us consider the action of
the product 
representation
\[
    h=(h_l,h_r):\pi_1(\Sigma)\rightarrow PSL(2,\mr)\times PSL(2,\mr)\,.
\]
on $\partial\mh^2\times\partial\mh^2$. Since neither $h_l$ nor $h_r$ fixes a point in
$\partial\mh^2$, $\pi_1(\Sigma)$ fixes no point on $\partial\mh^2\times\partial\mh^2$.

Given $\gamma\in\Gamma$ denote by $x_L^\pm(\gamma)$
(resp. $x_R^\pm(\gamma)$) the 
attractive and repulsive fixed points of $h_l(\gamma)$
(resp. $h_r(\gamma)$). If $h_l(\gamma)$ is parabolic, then 
 $x_L^+(\gamma)=x_L^-(\gamma)$ is the unique fixed point of $h_l(\gamma)$.
Let us introduce the following notations for the fixed points of $h(\gamma)$:
\[
   \begin{array}{ll} 
   p_{++}(\gamma)=(x_L^+(\gamma),x_R^+(\gamma)), &
    p_{+-}(\gamma)=(x_L^+(\gamma), 
    x_R^-(\gamma)),\\
    p_{-+}(\gamma)=(x_L^-(\gamma),x_R^+(\gamma)), &
    p_{--}(\gamma)=(x_L^-(\gamma), 
    x_R^-(\gamma))\,.
  \end{array}
\]

For every $\alpha,\gamma\in\pi_1(S)$ the following
identities hold 
\[
   \begin{array}{ll}
    p_{++}(\gamma^{-1})=p_{--}(\gamma), & p_{+-}(\gamma^{-1})=p_{-+}(\gamma),\\
    p_{\pm\pm}(\alpha\gamma\alpha^{-1})=h(\alpha)p_{\pm\pm}(\gamma)\,. &
   \end{array}
\]
It is also easy to see that for every
$p\in\partial\mh^2\times\partial\mh^2\setminus\{p_{-+}(\gamma),p_{+-}(\gamma), 
p_{--}(\gamma)\}$
\[
   \lim_{k\rightarrow+\infty} h(\gamma^k) (p)= p_{++}(\gamma)\,.
\]

A consequence of  the last fact is that any non-empty $h$-invariant closed subset of
$\partial\mh^2\times\partial\mh^2$ must contain $p_{++}(\gamma)$ for every
$\gamma\in\pi_1(\Sigma)$.

\begin{defi}
We define:
\[
   \Lambda=
   \Lambda(h_l,h_r)=\overline{\{p_{++}(\gamma)|\gamma\in\pi_1(\Sigma)\}}~. 
\] 
\end{defi}

By the remarks above, $\Lambda$ is the smallest non-empty closed $h$-invariant
subset of $\partial\mh^2\times\partial\mh^2$.

We state in the next proposition some basic properties of $\Lambda$, we
refer to \cite{mess,mess-notes,barbot-1,barbot-2} for the proofs. 

\begin{prop}\label{limitset:prop}
\begin{enumerate}
\item The projection on the first factor,
 $\pi_L:\partial\mh^2\times\partial\mh^2\rightarrow\partial\mh^2$, sends $\Lambda$ onto the
 limit set of $h_l$.
\item If both $\eta_l,\eta_r$ are complete then $\Lambda$ is the graph of the
homeomorphism of $f:\partial\mh^2\rightarrow\partial\mh^2$ conjugating $h_l$ and $h_r$.
Otherwise $\Lambda$ is a Cantor set. 
\item The metrics $\eta_l$ and $\eta_r$ are isotopic if and only if
  $\Lambda$ is contained in the boundary of a spacelike plane.
\end{enumerate}
\end{prop}
\begin{proof}
Notice that $\pi_L(\Lambda)$ is a closed subset of $\partial\mh^2$ invariant under $h_l$.
Thus it contains the limit set, say $\Lambda_l$, of $h_l$.

On the other hand, $\pi_L^{-1}(\Lambda_l)$ is a closed subset of
$\partial\mh^2$ invariant under $h$. So it must contain $\Lambda$.

Thus we have proved  
\begin{itemize}
\item $\Lambda_l\subset \pi(\Lambda)$
\item $\Lambda\subset \pi_L^{-1}(\Lambda_l)$
\end{itemize}
and we can conclude that $\Lambda_l=\pi_L(\Lambda)$.

For  point (2), let $f$ be the equivariant
homeomorphism of $\partial\mh^2$ conjugating $h_l$ with $h_r$.  Since the graph
of $f$ is invariant by $h$, it contains $\Lambda$. On the other hand,
since $\pi_L(\Lambda)=\Lambda_l=\partial\mh^2$, it follows that $\Lambda$
coincides with the whole graph.

Finally for the third point, notice that if $h_l$ and $h_r$ are
conjugated in $PSL_2(\mathbb R)$, then up to conjugation the points
$p_{++}(\gamma)$ lie all on the diagonal of $\partial\mh^2\times\partial\mh^2$.  It
follows that $\Lambda$ is contained in the boundary of $P_0$.

On the other hand if $\Lambda$ is contained in the boundary of some
spacelike plane $P_A$ then for every $\gamma\in\pi_1(S)$ we have
\[
      x^+_R(\gamma)= A x^+_L(\gamma)~.
\]
It follows that for every $\gamma$ the attractive and repulsive fixed points of
$h_l(\gamma)$ and $A^{-1}h_r(\gamma)A$ coincides. 
Thus we conclude that
\[
  h_l(\gamma)=A^{-1}h_r(\gamma)A~. 
 \]
\end{proof}

The  action of $\pi_1(\Sigma)$ on $\partial\mh^2\times\partial\mh^2$ is reminiscent of a
quasifuchsian action of $\pi_1(\Sigma)$ on 
$S^2=\partial_\infty \mh^3$. An important
difference with  that case is that the action of $\pi_1(\Sigma)$ on
$\partial\mh^2\times\partial\mh^2\setminus\Lambda$ is not proper. Indeed
$p_{+-}(\gamma)$  is not contained in $\Lambda$ and is fixed by $h(\gamma)$.
We are going to describe a maximal $h$-invariant domain of
$\partial\mh^2\times\partial\mh^2$, on which the action is properly discontinuous and causal.
This set could be regarded as the analogous of the discontinuity domain. 

Given an element $\gamma\in\pi_1(S)$ parallel to a puncture
we consider the two intervals $I_L(\gamma)$ and
$I_R(\gamma)$ in $\partial\mh^2$ that corresponds to
the infinite end of $\gamma$. If $h_l(\gamma)$ (resp. $h_r(\gamma)$) is
parabolic then $I_L(\gamma)$ (resp. $I_R(\gamma)$)
is reduced to a point.

\begin{prop}
The set $\Lambda$ is connectible by an achronal meridian.

The region, say $\mathcal G$, between the upper and lower meridians passing through
$\Lambda$ is
\[
     \overline{\bigcup_{\gamma\textrm{ parallel to a
           puncture}}I_L(\gamma)\times I_R(\gamma)}
\]
The action of $\pi_1(\Sigma)$on  $\mathring{\mathcal G}$ is free and properly discontinuous.
\end{prop}
\begin{proof}
Let $\Lambda_l^*$ be the set of conical limit points of $h_l$.  Let
$f:\mathcal H_l\rightarrow\mathcal H_r$ be the lifting of some
orientation-preserving diffeomorphism
$f:\Sigma_{\eta_l}\rightarrow\Sigma_{\eta_r}$.  We know that $f$
extends to a continuous map
\[
 f:\mathcal H_l\cup\Lambda_l^*\rightarrow\mathcal H_r\cup\partial\mh^2
\]
sending $\Lambda_l^*$ to some subset of $\Lambda_r$.

It is easy to see that the graph of $f|_{\Lambda_l^*}$ is contained in $\Lambda$.
We call this set $\Lambda^*$. Since it is invariant under $h$, its closure is
$\Lambda$.

So by Remark \ref{closure:rk} it is sufficient to prove that $\Lambda^*$ is connectible by
an achronal meridian.

Now let us take $x_1<x_2<x_3\in\Lambda_l^*$ and suppose $f(x_1)\neq
f(x_3)$.  Consider the oriented geodesic $l$ joining $x_1$ to $x_3$
and let $r$ be the half-line joining a point of $l$ to $x_2$.  By our
assumption, $r$ is contained on the right side bounded by $l$.

It follows that $f(r)$ is contained on the right side of $f(l)$. On the
other hand $f(l)$ has first end-point $f(x_1)$ and $f(r)$ joins a
point of $f(l)$ to $f(x_2)$. So $f(x_1)\leq f(x_2)\leq
f(x_3)$.

To conclude we have to show that $\mathcal G$ is
the closure of the union of $I_L(\gamma)\times I_R(\gamma)$

By remark \ref{rect:rk}, there is a one-to-one correspondence between rectangles of
$\mathcal G\setminus\Lambda$ and components of $\partial\mh^2\setminus\Lambda_l$.
Every component of $\partial\mh^2\setminus\Lambda_l$ is of form $I_L(\gamma)$ for
some peripheral element $\gamma$.

Let $G$ be the rectangle in $\mathcal G$ corresponding to the interval
$I_L(\gamma)\times I_R(\gamma)$. Since the interior of
$I_L(\gamma)\times I_R(\gamma)$ does not contain points of $\Lambda$,
by remark \ref{rect:rk} we have that $I_L(\gamma)\times I_R(\gamma)$
is contained in $G$.

On the other hand  notice that $p_{--}(\gamma)=(x_-^L(\gamma), x_-^R(\gamma))$ and
$p_{++}(\gamma)=(x_+^L(\gamma), x_+^R(\gamma))$ are both in
$\Lambda$. In particular this shows that if $C$ is a meridian curve
containing $\Lambda$ then $C\cap\pi_L^{-1}(I_L(\gamma))$ is contained
in $I_L(\gamma)\times I_R(\gamma)$.

This shows that $G\subset I_L(\gamma)\times I_R(\gamma)$
and the proof is complete.
\end{proof}

We says that a $\pi_1(\Sigma)$-invariant achronal meridian is
extremal if it is contained in the boundary of $\mathcal G$.

\begin{remark} \label{rk:number}
The number of non-degenerate rectangles $I_L(\gamma)\times I_R(\gamma)$ 
up to the action of $\pi_1(\Sigma)$ is equal to the number, say $k$, of punctures
of $\Sigma$ that corresponds to boundary components for both $\Sigma_{\eta_l}$ and
$\Sigma_{\eta_r}$.

It follows from remark \ref{mer:rm} that there are exactly $2^k$ extremal
$\pi_1(\Sigma)$-invariant  meridians.
\end{remark}

\section{Earthquakes and extremal invariant curves}

In this section we clarify the relation between earthquakes on a hyperbolic
surface with geodesic boundary and extremal
curves on the boundary at infinity $AdS_3$ which are invariant under the
action of a group. This will lead in particular to the proof of 
Theorem \ref{tm:earthquake}.

\subsection{Bent surfaces constructed from an earthquake}

Let us fix two admissible metrics $\eta_l$ and $\eta_r$.
Let
\[
 h=(h_l,h_r):\pi_1(\Sigma)\rightarrow PSL_2(\mr)\times PSL_2(\mr)~.
\]
be the representation whose components are the holonomies of
$\eta_l$ and $\eta_r$ respectively.

Let $\lambda$ be a measured geodesic lamination
on $\Sigma_{\eta_l}$ such that 
$E^r_\lambda(\Sigma_{\eta_l})=\Sigma_{\eta_r}$.
Consider the lifted earthquakes
$\tilde E:\tilde\Sigma_{\eta_l}\rightarrow\tilde{\Sigma}_{\eta_r}$ and  let
\[
b:\mathcal H\rightarrow AdS_3
\]
be the admissible $\pi_1(\Sigma)$-invariant past bent surface
constructed as in section \ref{earthtobend:sec}. 

We denote by $S_+$ the quotient of
$\mathcal H$ by the action of $\pi_1(\Sigma)$.
By Proposition \ref{pr:57}, $S_+=E^r_{\lambda/2}(\Sigma_{\eta_l})$.

\begin{prop} \label{pr:51}
The lower meridian passing through $\partial_\infty\mathcal H$, say $C$,
is an extremal $\pi_1(\Sigma)$-invariant meridian.
\end{prop}

\begin{proof}
Given an element $\gamma\in\pi_1(\Sigma)$ corresponding to a geodesic
boundary for both $\eta_l$ and $\eta_r$, let us set $G=I_L(\gamma)\times I_R(\gamma)$.
We have to show that $\partial_\infty\mathcal H$ does not intersect $\mathring G$.

Given a component $F$ of $\tilde\Sigma_\eta\setminus\lambda$ the corresponding stratum
of the bent surface $b(\mathcal H)$, say $K(F)$, is the ideal polygon in
$AdS_3$ whose end-points are the pairs $(x,Ax)$ where $x$ is an ideal point of $F$ and
$A\in PSL_2(\mathbb R)$ is determined by requiring that $\tilde E|_F=A|_F$.

Since $F$ is contained in the convex-hull of $\Lambda_l$, no ideal point of $F$ is contained
in $\mathring I_L(\gamma)$. It turns out that $K(F)\cap \mathring G=\emptyset$.
\end{proof}

Let us describe  describe  more precisely the curve $C$. In particular we will
describe for every region $G$ the intersection $G\cap C$ 

\begin{prop} \label{pr:82}
Let $\gamma\in\pi_1(\Sigma)$  
be the peripheral loop corresponding to a puncture $p\in\Sigma$.
Let $m$ be the total mass of $\lambda$ around $p$ 
and $a$ be the length of the
boundary component of $S_+$ corresponding to $\gamma$ ($a=0$ if $\gamma$
corresponds to a cusp in $S_+$).
Then $G\cap C$ is the lower curve if $m\leq a$, is the upper curve otherwise. 
\end{prop}

\begin{proof}
We can choose coordinates on $\mh^2$ -- considered in the Poincar\'e 
half-plane model -- in such a way that $h(\gamma)z= e^{2a}z$
($h(\gamma)z=z+1$ if $a=0$). In particular if $a\neq 0$ we can suppose that
$\tilde F\subset\{z\in\mc| Re(z)<0\}$.


If $m\neq 0$, we have two possibilities: leaves near $p$ lift to leaves with
endpoints either at $0$ or at $\infty$ (if $a=0$ only the last possibility
holds, in the other cases the choice depends on the way of spiraling of
$\lambda$ around $p$).
A hyperbolic transformation with attractive fixed point at $\infty$
(resp. $0$) is upper triangular (resp. lower triangular).

Thus if we choose a base-point near the puncture, it is easy to see that if
leaves of $\lambda$ near $p$ lift to geodesics with an end-point at $\infty$, then 
as in Proposition~\ref{pr:length} we have

\begin{eqnarray}\label{holonomies}
\begin{array}{cc}
h_l(\gamma)=\left(\begin{array}{cc} e^{a-m} & *\\ 0 &
    e^{-(a-m)}\end{array}\right)
&
h_r(\gamma)=\left(\begin{array}{cc} e^{a+m} & *\\ 0 &
    e^{-(a+m)}\end{array}\right)
\end{array}
\end{eqnarray}

In the same way, if the common endpoint of these geodesics is $0$ then
we have
\[
\begin{array}{cc}
h_l(\gamma)=\left(\begin{array}{cc} e^{(a+m)} & 0\\ * &
    e^{-(a+m)}\end{array}\right)
&
h_r(\gamma)=\left(\begin{array}{cc} e^{a-m} & 0\\ * &
    e^{-(a-m)}\end{array}\right)\,.
\end{array}
\]
Since $h_l(\gamma)$ and $h_r(\gamma)$ are assumed to be 
hyperbolic, we find that $a\neq m$.

Let us distinguish three cases.
\begin{enumerate}
\item 
$m=0$: in this case $\gamma$ corresponds to a boundary
   component of $F$ and  the bending map 
\[
  \beta:\tilde F\rightarrow AdS_3
\]
extends  to the axis of $\gamma$, say $l$. Moreover, the image of $l$ is the
axis of $(\gamma)$ (that is the geodesic joining the limit end-points of
$G$). Thus $\partial F_\lambda\cap G$ is the past extremal curve.

\item $0<m<a$. Let $l_0$ be a leaf of $\tilde \lambda$ with an endpoint fixed by
   $\gamma$, say $s_0\in\{0,\infty\}$. Denote by $t_0\in\mr_{<0}$ its other
   endpoint. Notice that the image of $\gamma^k(l_0)$ through $\beta$ is a
   geodesic, say $\hat l_k$, of $AdS_3$, with end-points equal to $(s_0, s_0)$
   and to $(h_l^k(\gamma) t_0, h_r^k(\gamma)t_0)$.
   
   Notice that  $(h_l^k(\gamma) t_0, h_r^k(\gamma)t_0)$ converges to
   $p_{++}(\gamma)$ as $k\rightarrow +\infty$ and to $p_{--}(\gamma)$ as
   $k\rightarrow-\infty$. On the other hand, since $m<a$ we have that
   $(s_0,s_0)=p_{++}(\gamma)$ if $s_0=\infty$ and $(s_0,s_0)=p_{--}(\gamma)$
   if $s_0=0$.

   Suppose for the sake of simplicity that $s_0=\infty$. Then the
   geodesic $\hat l_k$ converges to $(s_0,s_0)$ as $k\rightarrow+\infty$, but
   $\hat l_k$ converges to the geodesic $c_\gamma$ with end-points $p_{--}(\gamma)$ and
   $p_{++}(\gamma)$ as $k\rightarrow -\infty$. It follows that such a geodesic
   is contained in the boundary of $\mathcal H$. Then, as in the
   previous case, we get that $\partial F_\lambda\cap G$ is the past extremal
   curve.

\item $m>a$. In this case $(s_0,s_0)\in\{p_{-+}(\gamma),
  p_{+-}(\gamma)\}$. Thus, by arguing as before, we get that $\hat
  l_k$ converges to the lightlike segment in $\partial_\infty AdS_3$
  joining $p_{++}(\gamma)$ to $(s_0,s_0)$ as $k\rightarrow+\infty$ and
  converges to the lightlike segment in $\partial_\infty AdS_3$
  joining $p_{--}(\gamma)$ to $(s_0,s_0)$ as $k\rightarrow -\infty$.
  Thus the upper extremal curve of $G(\gamma)$ is contained in the
  closure of the image of $\beta$ and thus it is contained in $\tilde
  F_\lambda$.

\end{enumerate}
Thus we have proved that if $(F,\lambda)$ encodes a bent surface in $M$, and $F$
is admissible, then the boundary curve of $\tilde F_\lambda$ is extremal.
\end{proof}

\subsection{Constructing an earthquake from an invariant curve}

In this subsection we prove the following proposition,
the last missing tool for the proof of Theorem \ref{tm:earthquake}.

\begin{prop}\label{prop:51bis}
Let $C$ be an extremal $\pi_1(\Sigma)$-invariant achronal meridian
and let $K$ be its convex hull.
There is an earthquake $E^r_\lambda$ sending $\Sigma_{\eta_l}$ into
$\Sigma_{\eta_r}$ whose associated admissible bent surface
is the spacelike part of $\partial_+K$.
\end{prop}


To prove Proposition \ref{prop:51bis} we need the following simple
technical lemma of hyperbolic geometry.

\begin{lem}\label{conn:lm}
Let $H$ be a  closed path-connected subset of  $\mathbb H^2$ and 
suppose that there exists a decomposition:
\[
   H=\bigcup_i H_i
\]
where each $H_i$ is either an ideal geodesic polygon or a geodesic,
and $H_i$ and $H_j$ are weakly separated in $\mathbb H^2$.
Then $H$ is convex.
\end{lem}

\begin{proof}
Let $c:[0,1]\rightarrow\mh^2$ be the geodesic segment joining two
points $p,q\in H$.  Consider the set of $t$ such that $c|_{[0,t]}$
is contained in $H$ and let $t_0$ the $\sup$ of this set.

By contradiction suppose $t_0<1$ and let $p'=c(t_0)$.  It is not
difficult to construct a family of geodesics $l_n$ such that each
$l_n$ lies in the frontier of some $H_i$ and $l_n$ intersects $c$ at a
point $c(t_n)$ with $t_n\rightarrow t_0^-$.  Up to a subsequence,
$l_n$ converges to a geodesic $l$ that has the property that is weakly
separated from all $H_i$.
Moreover $p$ and $q$ lie on opposite sides of $l$.
Let $U$ be the open half-plane bounded by $l$ and containing $q$.

Thus if $\gamma$ is the path joining $p$ to $q$ in $H$, then there is
a time $s_0<1$ such that $\gamma(s_0)\in l$ and $\gamma(s)$ lies in
$U$ for $s\in(s_0,1]$. Let us consider the set $L$ of points
$s\in[s_0,1]$ such that $\gamma(s)$ lies in the frontier of some
$H_i$.  Notice that $L$ is closed
and that the components of $(s_0,1]\setminus L$ are contained in the interior.  
Now if the set $L$ does not
accumulate on $s_0$, this means that there exists $H_i$ containing all
the path $\gamma|_{[s_0,s_0+\delta]}$. Since $H_i$ is weakly separated by $l$, it turns out that
$H_i$ is an ideal geodesic polygon and $l$ is one of its boundary component. But then
$c(t)\in H_i$ for $t\in[t_0, t_0+\delta']$, contradicting the definition of $t_0$.

Suppose now that $s_0$ is an accumulation point for $L$. There is a
$\delta>0$ such that if $s\leq s_0+\delta$ and $s\in L$ then
$\gamma(s)$ lies in the frontier $l(s)$ of some $H_i$ that must
intersect $c$ at some point $c(t)$ with $t>t_0$, and $t\rightarrow
t_0^+$ as $s\rightarrow s_0$. We can also suppose that $s_0+\delta\in L$
and let us set $\delta'$ such that $c(t_0+\delta')\in l(s_0+\delta)$.
 
The set of points 
$\hat L=\{t|c(t)\in l(s)\}$ is a closed subset that accumulates at $s_0$.
On the other hand if $I$ is a component of $[t_0,t_0+\delta]\setminus L$
then its end-points lies respectively in $l(s)$ and $l(s')$ and
$(s,s')$ is a component of $[s_0,s_0+\delta]\setminus L$.
It follows that $l(s)$ and $l(s')$ are boundary components of some $H_i$, but then
the segment $c(I)$ is contained in $H_i$.
This shows that $c|_{[t_0,t_0+\delta']}$ is contained in $H$ and we get a contradiction.
\end{proof}


\begin{proof}[Proof of Proposition \ref{prop:51bis}]
Let us consider the set
\[
  \mathcal P=\{A\in PSL_2(\mr)| P_A\textrm{ is a support plane for }\partial_+K\}~.
\]
If $A\in\mathcal P$ then 
$P_A\cap\partial_+K$ is the convex hull of some set
$\{x_i,Ax_i\}_{i\in I}=C\cap P_A$.

We define the stratum associated to $A$ to be the convex hull, say $F(A)$, 
in $\mh^2$ of $\{x_i\}_{i\in I}$.
By Remark \ref{intersection:rk}, $F(A)$ is weakly separated from $F(A')$.

By the invariance of $\partial_+K$,
if $P$ is a support plane for $\partial_+K$, then so is 
$(h_l(\gamma),h_r(\gamma))(P)$. Thus if $A$ lies  in $\mathcal P$ then
so does $\gamma\cdot A=h_r(\gamma)\circ A\circ h_l(\gamma)^{-1}$ and
$F(\gamma\cdot A)= h_l(\gamma)(F(A))$.
In particular, the set
\[
\mathcal H=\overline{\bigcup_{A\in\mathcal P} F(A)}
\]
is $h_l$-invariant. 

We claim that $\mathcal H$ is the convex core of
$h_l$.  First we prove that $F(A)$ is contained in the convex core for
every $A\in\mathcal P$.  By contradiction, if $F(A)$ is not contained
in the convex core there is some ideal point of $F(A)$ that is not
contained in $\Lambda_l$.  Thus there is some point $(x,y)\in P_A\cap
C$ such that $x\in \mathring{I_l(\gamma)}$ for some peripheral
$\gamma$ Since $C$ is an extremal meridian the point $(x,y)$ must lie
on the interior of some edge $e$ of $I_l(\gamma)\times I_r(\gamma)$.
Since $P_A$ is a support plane for $K$, it cannot disconnect $C$ (more
precisely, in a suitable affine chart $C$ is contained in one of the
two closed half-spaces bounded by $P_A$).  It follows that $P_A$ must
contain $e$. Since $P_A$ is spacelike, this gives a contradiction
(spacelike planes intersect every leaf only once).

This proves that $\mathcal H$ is contained in the convex hull.
To prove the reverse inclusion, we only have to check  that 
$\mathcal H$ is convex.

By Lemma \ref{conn:lm} it is sufficient to prove that
$\mathcal H$ is path-connected.  Notice that by definition
$r_A(F_A)=P_A\cap\partial_+K$.  So given two points $p\in F(A)$ and
$q\in F(A')$ let us consider the corresponding points
$\hat p=r_A(p)$ and $\hat q=r_A(q)$ in $\partial_+K$.
By classical facts on convex subset in $\mathbb R^3$, 
there is a continuous path
\[
   u:[0,1]\rightarrow T(AdS_3),\qquad u(t)=(\hat p(t),v(t))
\]
such that $\hat p(t)$ is a path in the spacelike part of 
$\partial_+K$ joining $\hat p$ to $\hat q$
and $v(t)$ is a vector orthogonal to some support plane $P(t)$ at 
$\hat p(t)$ of $\partial_+K$.

Let $A(t)\in PSL_2(\mathbb R)$ be such that $P(t)=P_{A(t)}$.
Then the path
\[
    p(t)= r_{A(t)}^{-1}p(t)
\]
is a continuous path in $\mathbb H^2$ joining $p$ to $q$.  Since
$p(t)\in F(A(t))$ we conclude that $p$ and $q$ are connected by an arc
in $\mathcal H$, so it is connected.

We  consider on $\mathcal H$ the geodesic lamination
\[
   L=\bigcup_{A\in\mathcal P: F(A)\textrm{is a geodesic}} F(A)\quad\cup\quad
  \bigcup_{A\in\mathcal P: F(A)\textrm{is a polygon}}\partial F(A)
\] 

We construct a right earthquake on $\mathcal H$  with fault locus $L$
and such that the corresponding bent surface is $\partial_+K$.

Indeed every stratum $F$ of $L$ coincides with $F(A)$ for some
$A\in\mathcal P$. So we can select for every stratum $F$ an element
$A=A(F)\in\mathcal P$ such that $F=F(A)$ 
(the choice is unique almost everywhere).
So we define the map
\[
   E:\mathcal H\rightarrow\mathbb H^2
\]
such that $E|_{F}=A(F)$

Let us consider two strata $F$ and $F'$. 
The planes $P_{A(F)}$ and $P_{A(F')}$ meet along a line $l$.
Orient $l$ in such a way that the signed angle between $P_{A(F)}$ and
$P_{A(F')}$ is positive.
Let $(x_-,y_-)$ and $(x_+,y_+)$ be the end-points of $l$.
Since $y_\pm=A(F)x_\pm=A(F')x_{\pm}$ it turns out that
$x_-$ and $x_+$ are fixed points for the comparison isometry
$B^*=A(F)^{-1}\circ A(F)$.

Thus $B^*$ is a hyperbolic translation whose axis is the geodesic $s$ in
$\mathbb H^2$ with end-points $x_-$ and $x_+$.  From Remark \ref{intersection:rk}
$s$ separates $F$ from $F'$.  So in order to conclude we just
have to check that $F'$ is moved by $B^*$ on the right as viewed from
$F$. Since $F'$ is contained in the left side bounded by $s$, it is
sufficient to prove that $B^*$ acts by a positive translation on $s$.

Notice that the isometry $(B^*,1)$ preserves  $P_{A(F')}$ into $P_{A(F)}$.
Since the signed angle $\alpha(P_{A(F')},P_{A(F)})<0$, then the rotation component of
$(B^*,1)$ is negative.
 By Lemma \ref{comp:lm} we conclude that $t(B^*)>0$.

By applying construction of Proposition \ref{pr:618} to the earthquake $E$, we
easily check
 that the bent surface associated to $E$ is the spacelike part of $\partial_+K$.
\end{proof}

\subsection{Proof of the main result}

We are now ready to prove the earthquake theorem for hyperbolic surfaces
with geodesic boundary.

\begin{proof}[Proof of Theorem \ref{tm:earthquake}]
Proposition \ref{prop:51bis}, along with Proposition \ref{pr:51}, shows
that there is a 1-to-1 correspondence between right
earthquakes relating $\eta_l$ to $\eta_r$
and extremal curves. 
According to Remark \ref{rk:number}, the number of extremal curves is equal to 
$2^k$, where $k$ is the number of boundary components which are cusps
neither for $\eta_l$ nor for $\eta_r$. 
Theorem \ref{tm:earthquake} follows.
\end{proof}

\section{The enhanced Teichm\"uller space}

The enhanced Teichm\"uller space $\cTh_{g,n}$ of a compact surface 
with boundary 
is defined in the introduction (Definition \ref{df:enhanced}). 
There is a natural topology on $\cTh_{g,n}$, which restricts to
the domains where all $\epsilon_i$ are constant as the usual 
topology on $\cT_{g,n}$. It can be defined through a family of 
neighborhoods of a point, involving quasiconformal homeomorphisms
diffeomorphic to the identity, as well of course as the $\epsilon_i$ 
(we leave the details to the reader).

We now
wish to define earthquakes as maps from $\cTh_{g,n}$ to itself. A naive
possibility would be to define it as the earthquakes on $\cT_{g,n}$, 
adding some information on the signs assigned to boundary components.
This however would yield a definition which is not quite satisfactory,
since right earthquakes would not have some desirable properties, like
those appearing in Proposition \ref{pr:flow} below.

\subsubsection*{Reflections of geodesic laminations}

Some preliminary definitions are needed. Here we consider a compact
surface $S$ of genus $g$ with $n$ boundary components, and a hyperbolic
metric $h$ with geodesic boundary on $S$. 
Let $c_0$ be one of the boundary components of $S$ which is not a cusp.
Let $\gamma$ be a complete oriented embedded geodesic in $(S,h)$ which is 
asymptotic to $c_0$ on its positive endpoint, i.e., which spirals onto $c_0$
as $t\rightarrow \infty$.

\begin{defi}
The {\bf reflection} of $\gamma$ relative to $c_0$, denoted by 
$\sigma_{c_0}(\gamma)$, is the geodesic 
in $(S,h)$ obtained as follow. Let $\gammab$ be any lift of $\gamma$
to the universal cover $\St$ of $S$, so that $\gammab$ has its endpoint on 
the positive side at an endpoint of a lift $\cb_0$ of $c_0$. We define 
$\sigma_{c_0}(\gamma)$ to be the projection on $S$ of the complete geodesic having as 
its positive endpoints
the other endpoint of $\gammab$ and the other endpoint of $\cb_0$.
\end{defi}

Note that considering an oriented geodesic here is necessary only if $\gamma$
spirals on $c_0$ at both ends. The existence of this reflected geodesic 
can also be considered in light of Lemma \ref{b:lem}.

\begin{remarks}
\begin{enumerate}
\item $\rho_{c_0}(\gamma)$ is also embedded,
\item if $\gamma_1$ and $\gamma_2$ are two geodesics asymptotic to 
$c_0$ which are disjoint, then 
$\rho_{c_0}(\gamma_1)$ and $\rho_{c_0}(\gamma_2)$ are also disjoint,
\item if $\gamma_1$ is a geodesic asymptotic to $c_0$ and $\gamma_2$
is a geodesic not asymptotic to $c_0$, and if $\gamma_1$ and $\gamma_2$
are disjoint, then $\sigma_{c_0}(\gamma_1)$ and $\gamma_2$ are disjoint.
\end{enumerate}
\end{remarks}

\pf
For the first point let $(\gamma_t)_{t\in (0,1)}$ be a one-parameter
family of geodesic rays starting from a point of $c_0$ and ending
at the common endpoint of $\gamma$ and of $\sigma_{c_0}(\gamma)$,
such that $\lim_{t\rightarrow 0}\gamma_t=\gamma$, 
$\lim_{t\rightarrow 1}\gamma_t =\sigma_{c_0}(\gamma)$.
Since $\gamma_0=\gamma$ is embedded, it is not difficult to show
that $\gamma_t$ is embedded for $t$ small enough. Suppose that 
$\gamma_t$ is not embedded for some $t\in (0,1)$, and let $t_0$
be the infimum of the $t\in (0,1)$ such that $\gamma_t$
is not embedded, then $\gamma_{t_0}$ would have a self-tangency
point, which is impossible. So $\gamma_t$ is embedded for all
$t\in (0,1)$, and therefore $\gamma_1=\sigma_{c_0}(\gamma)$ is also
embedded. This proves the first point.

For the second point let $\gammab_1$ and $\gammab_2$ be any lifts
of $\gamma_1$ and of $\gamma_2$ to $\St$. Since $\gamma_1$ and 
$\gamma_2$ are disjoint, their lifts $\gammab_1$ and $\gammab_2$
are also disjoint. $\sigma_{c_0}(\gamma_1)$ is the image by the
projection $\St\rightarrow S$ of the geodesic, which we can denote
by $\sigma_{c_0}(\gamma_1)$, which has one endpoint in common with
$\gammab_1$ (not on a lift of $c_0$) while the other is an endpoint
of a lift of $c_0$ which has as its other endpoint an endpoint of 
$\gammab_1$. The same description holds for $\gammab_2$. It follows
that $\sigma_{c_0}(\gammab_1)$ and $\sigma_{c_0}(\gammab_2)$ are 
also disjoint. Since this is true for all lifts of $\gamma_1$
and $\gamma_2$ to $\St$, $\sigma_{c_0}(\gamma_1)$ and 
$\sigma_{c_0}(\gamma_2)$ are disjoint, which proves the second point.

The third point can be proved using the same argument, we leave 
the details to the reader.
\cvd

\begin{defi}
Let $\lambda$ be a measured geodesic lamination on $(S,h)$. The
{\it reflection} of $\lambda$ is the measured lamination, denoted 
by $\sigma_{c_0}(\lambda)$, obtained
by replacing each leave of $\lambda$ which is asymptotic to 
$c_0$ by its reflection relative to $c_0$.
\end{defi}

The previous remarks makes this definition possible, since they
show that $\sigma_{c_0}(\lambda)$ is again a measured geodesic lamination
(its support is a disjoint union of geodesics). Note that the reflection
map, acting on measured geodesic laminations, has some simple properties:
\begin{itemize}
\item if $\lambda$ is any measured geodesic laminations on $S$
then $\sigma_{c_0}^{2}(\lambda)=\lambda$,
\item if $c_0$ and $c_1$ are two boundary components of $S$ 
then $\sigma_{c_0}$ and $\sigma_{c_1}$ commute.
\end{itemize}
In both cases the proofs follow by considering the corresponding
statements for geodesics.

\subsubsection*{Earthquakes on the enhanced Teichm\"uller space}

The definition of earthquakes on $\cTh_{g,n}$ is based on the
reflection of measured geodesic laminations.

First we define the earthquake on the subset of $\cTh_{g,n}$, say $\cTh_{g,n}'$ of admissible
metrics without cusps.
\begin{defi} \label{df:e_r}
The map $E_r:\cML_{g,n}\times \cTh_{g,n}'\rightarrow \cTh_{g,n}$ is
defined as follows.
Let $(\eta, \epsilon_1, \cdots, \epsilon_n)\in \cTh_{g,n}'$, and let
$\lambda\in \cML_{g,n}$. Consider the measured lamination 
$\lambdab$ obtained by taking the reflection of $\lambda$
with respect to all boundary components $c_i$ of $S$ for which 
$\epsilon_i=-1$, and let $\bar \eta = E_r(\lambdab)(\eta)$. Finally, for
$i=1,\cdots,n$, let $\epsilonb_i=\epsilon_i$ if the right earthquake
$E_r(\lambdab)$ does not change the direction in which $\lambdab$
spirals into $c_i$, $\epsilonb_i=0$ if $c_i$ is become a cusp 
and $\epsilonb_i=-\epsilon_i$ otherwise. Then
$E_r(\lambda)(\eta, \epsilon_1, \cdots, \epsilon_n) = 
(\bar \eta, \epsilonb_1, \cdots, \epsilonb_n)$.
\end{defi}

The following lemma ensures that it is possible to extend $E_r$ to the whole
$\cTh_{g,n}$.
\begin{lem}
For any $\lambda\in\Mm\Ll_{g,n}$ the map
$E_r(\lambda):\cTh_{g,n}'\rightarrow\cTh_{g,n}$ 
extends continuously on $\cTh_{g,n}$.
\end{lem}

\pf
Given a point $(\eta,\epsilon_1,\ldots,\epsilon_k)$ corresponding to a metric 
$\eta$ with $n-k$ cusps, denote by $\lambda'$ the measured
geodesic lamination of $\Mm\Ll_{g,n}(\eta)$ corresponding to $\eta$. 
Let $\bar\lambda'$ be the lamination obtained by reflecting $\lambda'$ with
respect to all boundary components for which $\epsilon_i=-1$ and let $\bar
\eta=E_r(\lambdab')(\eta)$. For $i=1,\ldots, k$ the sign $\bar\epsilon_i$ is
defined as in the previous case, whereas for $i=k+1,\ldots, n$ the sign
$\epsilonb_i=1$ if the lamination $\lambda$ (that is a lamination of a surface
without cusp) spiral in the positive way with respect to the standard
spiraling orientation, and $\epsilon_i=-1$ otherwise.
Finally let us define
\[
   E_r(\lambda)(\eta)=(\bar\eta,\bar\epsilon_1,\ldots,\bar\epsilon_n)\,.
\]

It is clear that the composition 
\begin{equation}\label{co:eq}
\begin{CD}
   \cTh_{g,n}@> E_r(\lambda)>>\cTh_{g,n}@>\pi>>\cT_{g,n}
\end{CD}
\end{equation}
is continuous.
To conclude it is then sufficient to show the following points:
\begin{enumerate}
\item 
 if $c_i$ is a cusp with respect to $\eta_0$, then there is a neighbourhood of
  $\eta$ such that $\bar\epsilon_i$ is constant.
\item
 if $c_i$ is a cusp with respect to $E_r(h)(\eta_0)$, then for $a>0$ there is a
  neighbourhood of $\eta_0$ such that the length of $c_i$ with respect 
  $E_r(\hat\lambda)(\eta)$ is smaller than $a$ for $\eta$ in that neighborhood.
\end{enumerate}

The second point follows from the continuity of (\ref{co:eq}).
For the first point, let $\eps=1$ if $\lambda$ spirals in the positive way around
$c_i$ and $\eps=-1$ otherwise. Notice that in a small neighbourhood $U$ of
$\eta_0$ (precisely the set of $\eta$ for which the length of $c_i$ is less
than the $\lambda$ mass of $c_i$) the lamination $E_r(\lambda)$ spirals in the positive way around $c_i$
Thus if you take a point $(\eta,\epsilon_1,\ldots,\epsilon_n)\in U$ then the
corresponding $\overline\epsilon_i$ is equal to $\eps$.
\cvd

There is a corresponding definition of left earthquake, and we
call $E_l:\cML_{g,n}\times \cTh_{g,n}\rightarrow \cTh_{g,n}$. It
follows directly from the definition that, for all $\lambda\in
\cML_{g,n}$, $E_l(\lambda)=E_r(\lambda)^{-1}$. 

As already mentioned above, Definition \ref{df:e_r} has some
desirable properties, that can not easily be achieved by more
simple-looking definitions. 

\begin{prop} \label{pr:cont}
The map $E_r:\cML_{g,n}\times \cTh_{g,n}\rightarrow \cTh_{g,n}$ is
continuous.
\end{prop}

\pf
We have seen above that the variation of the length of a boundary
component under $E_r(\lambda)$ is proportional to the mass of 
$\lambda$. The proof is therefore a direct consequence of the 
definition of the topology on $\cTh_{g,n}$.
\cvd

\begin{prop} \label{pr:flow}
Let $\lambda\in \cML_{g,n}$, and let $t,t'\in \R_{>0}$.
Then $E_r((t+t')\lambda) = E_r(t'\lambda)\circ E_r(t\lambda)$.
\end{prop}

\pf

Let $(h, \epsilon_1, \cdots, \epsilon_n)\in \cTh_{g,n}$, and let 
$(h', \epsilon'_1, \cdots, \epsilon'_n) =
E_r(t\lambda)(h, \epsilon_1, \cdots, \epsilon_n)$. 
Let  $\lambdab$ be the image of $\lambda$ under
the reflection relative to all boundary components of $S$ for which
$\epsilon_i=-1$. 
Suppose first that $\forall i\in \{ 1,\cdots, n\}, \epsilon'_i=\epsilon_i$.
This implies that, after the right earthquake $E_r(t\lambda)$
(considered as an earthquake acting on $\cTh_{g,n=}$)
$\lambdab$ spirals in the same direction onto each of the boundary
components of $S$. In other terms $\lambdab$ remains the
same measured lamination after the earthquake $E_r(\lambda)$. The fact that 
$$ E_r(t'\lambda)\circ E_r(t\lambda)(h, \epsilon_1, \cdots, \epsilon_n) 
= E_r((t+t')\lambda)(h, \epsilon_1, \cdots, \epsilon_n) $$
then follows directly from the definition of an earthquake (through
the right-quake cocycle as seen in section 3).

Suppose now that some of the $\epsilonb_i$ are different from the 
corresponding
$\epsilon_i$, and let $\lambdab'$ be the image of $\lambda$ under
the reflection relative to all boundary components of $S$ for which
$\epsilonb_i=-1$. The definition of the image of an element of 
$\cTh_{g,n}$ implies that $\lambdab'$ is the image of $\lambdab$
under the earthquake $E_r(t\lambda)$ -- the boundary components
of $S$ for which $\lambdab$ and $\lambdab'$ circle in
opposite directions are precisely those for which $\epsilonb_i
\neq \epsilon_i$. Thus the result follows again from an elementary
argument based on the right-quake cocycle.
\cvd

\subsubsection*{Proof of Theorem \ref{tm:earthquake2}}

The results of section 8 show that there are $2^k$ left earthquakes sending
a given hyperbolic metric (considered as an element of $\cT_{g,n})$ to another
one, where $k$ is the number of punctures corresponding to geodesic
boundary components in both metrics. The corresponding measured laminations
are the bending lamination of the future boundary of a convex retract $U$
of $M$ which has as boundary curve an extremal curve (see Proposition
\ref{pr:51}). It was also noted (in Proposition \ref{pr:82}) 
that each boundary curve is the upper boundary curve when the (signed)
mass $m$ of the measured lamination at the corresponding boundary component
of $F$ is bigger than $a$, the length of that boundary component in the induced
metric on the bent surface in $M$. However we have seen in section 3
that $m>a$ if and only if the lamination $\lambda$ spirals in the opposite 
direction on that boundary component after the left earthquake is performed.
So each upper extremal curve corresponds to a boundary component for which the
spiraling orientation is reversed by the earthquake along $\lambda$,
while each lower extremal curve corresponds to a boundary component for 
which the spiraling orientation remains the same. This proves Theorem 
\ref{tm:earthquake2}.

\section{Multi Black Holes}

By now we have studied the action of $h$ on the boundary of $AdS_3$. Let us
now consider the action of $h$ on $AdS_3$. A first easy remark is that such an
action is neither proper nor causal. For instance, the lightlike plane $P$ 
which is tangent to the boundary at infinity of $AdS_3$ at $p_{++}(\gamma)$ 
(considered in the projective model of $AdS_3$) is preserved
by $h(\gamma)$, and the orbits of $h(\gamma)$ on $P$ are contained in lightlike
rays. Moreover notice that if $h_l=h_r$ then $id$ -- considered as an element of $AdS_3$,
identified with $PSL(2,\R)$ as explained above -- is a fixed point for $h$.

In \cite{barbot-1,barbot-2} 
it has been shown that there exists a maximal domain, say
$\Omega=\Omega(h_l,h_r)$ of $AdS_3$ such that the action of $h$ on
$\Omega$ is free and properly discontinuous and the quotient
$\Omega/h(\pi_1(\Sigma))$ is a strongly causal Lorentzian manifold
homeomorphic to $\Sigma\times\mr$. We will refer to this quotient
as a MBH spacetime $M=M(h)=\Omega/h(\pi_1(\Sigma))$.

Let $\Kk$ be the convex hull of the limit set $\Lambda$ in $AdS_3$.
Recall that given an oriented space-like plane $P$ in $AdS^3$, all
time-like geodesic planes orthogonal to $P$ intersect at distance
exactly $\pi/2$ from $P$. The intersection point on the positive
side of $P$ is called the point dual to $P$.   
Using this notion, $\Omega$ can be defined as the set of points whose dual planes 
are disjoint from $\Kk$ (see \cite{barbot-1,barbot-2} for details).

Let us collect some properties of $\Omega$:
\begin{enumerate}
\item It is convex and strongly causal.
\item The intersection with the boundary at infinity of $AdS_3$ of the
  closure of $\Omega$ is the asymptotic region 
  $\Gg$ of $h$, as described in the previous subsection.
\item It contains the convex core of $\Kk$.
\end{enumerate}
Notice that $\Omega$ is not globally hyperbolic. Any globally hyperbolic 
spacetime with holonomy equal to $h$ isometrically embeds into $\Omega$.
Thus, $\Omega$ can also be described as the union of all $h$-invariant  
globally hyperbolic domains. Such domains, in turn, are in one-to-one
correspondence with $h$-invariant no-where timelike closed curves in 
$\partial AdS_3$. Let us note that $\Kk$ is contained in the convex hull of every 
such closed curve in $\partial AdS_3$.

The ``black hole'' of $\Omega$ is, by definition, the set of points
that cannot be connected to $\Gg$ along any future-directed causal path 
(i.e. the domain in $\Omega$ that is
causally disconnected from the ``infinity'' $\Gg$ in the future). 
Barbot \cite{barbot-1,barbot-2}
pointed out that this set is globally hyperbolic and corresponds to the 
extremal curve in the boundary obtained by choosing the arc $\wedge$ in each 
$\Gg(\gamma)$.

There is also a ``white hole'', that is the set of points that cannot be
connected to $\Gg$ along any \emph{past-directed} causal curve.
It is the globally hyperbolic domain whose boundary at infinity
is obtained by choosing the extremal arc $\vee$ in each $\Gg(\gamma)$.

The intersection of the black hole and the white hole is the set of points
disconnected from $\Gg$ both in the future and in the past. It can be regarded
as the set of points contained in all $h$-invariant MGH domains (in particular
it contains $\Kk$).


Notice that for each $\gamma\in\pi_1(\Sigma)$ corresponding to a
 non-degenerate AR, the
geodesic, say $c_\gamma$, joining $p_{--}(\gamma)$ to $p_{++}(\gamma)$ in
 $AdS_3$ is contained in 
the boundary of $\Kk$. This geodesic is contained in the lightlike planes dual
to $p_{+-}(\gamma)$ and $p_{-+}(\gamma)$. Consider then the lightlike
triangles with base $c_\gamma$ and vertex respectively in $p_{-+}(\gamma)$ and
$p_{+-}(\gamma)$. The union of these triangles disconnects $\Omega$ in two
 regions, the one that faces $\Gg(\gamma)$ is called the asymptotic region of
 $\gamma$ in $\Omega$.

The union of such triangles disconnects $\Omega$ in an
``internal'' piece, that is $\hat\Omega$, and a certain set of regions that
 faces the non-degenerate AR's. We call such regions the asymptotic regions of
 $\Omega$. The asymptotic regions of the MBH spacetime are defined
as the corresponding quotients.

The domain $\hat\Omega$ turns out to be the union of the black hole and the white
hole.
Thus the boundary of $\hat M=\hat\Omega/h$ is formed by $k$ annuli, each of
which is the union of two lightlike totally geodesic annuli along a spacelike
geodesics. Notice that $\hat M$ is a strong deformation retract of $M$.

The set of geodesics $\{c_{\gamma}\}$ is contained in the boundary of
$\partial\Kk$.
They disconnect the boundary in two bent surfaces, the upper and the lower
boundary of $\partial\Kk$.  The intrinsic metric on them is
hyperbolic, and in fact they are isometric to some straight convex sets of
$\mh^2$. Moreover the bending gives rise to a measured geodesic lamination
on each.
Clearly $\partial_\pm \Kk$ are invariant under the action of $h$, and
$\partial_\pm\Kk/h$ produces an admissible structure.

\section{Some remarks}

\subsubsection*{Ends versus cone singularities.}

The statements presented here, concerning earthquakes on hyperbolic surfaces 
with geodesic boundary, can quite naturally be compared to corresponding results
on closed surfaces endowed with hyperbolic metrics with cone singularities (as
in \cite{cone}). Indeed cone singularities can in a fairly natural way be 
considered as analytic continuations of geodesic boundary components when 
the length becomes imaginary. Another way to state the relation between 
the two is that black holes (or more precisely, singularities inside them, i.e.
bending lines on the boundary of $M$) are ``particles'' moving along
spacelike geodesics.

However this analogy has limits. One of them is that the ``earthquake theorem'' 
of \cite{cone} keeps the angle at the cone singularities fixed, so that two 
metrics are related by a unique right earthquake and there is 
no analog of the appearance of the enhanced Teichm\"uller space, which is a
key feature for hyperbolic surfaces with geodesic boundary.

There might very well be a statement generalizing both the main result here
and the main result of \cite{cone}, and describing the earthquakes between
two hyperbolic metrics having both cone singularities (of angle less than
$\pi$) and geodesic boundary components. One could even imagine a proof
based on a Mess type parametrization, by a right and left hyperbolic metric,
of the space of multi-black holes of a given topology containing ``particles''
of fixed cone angle. 

\subsubsection*{Other possible proofs.}

There are at least two possible proofs of Thurston's Earthquake Theorem for
(smooth) hyperbolic metrics on closed surfaces, in addition to the Mess
argument used above. One, originating in work of Kerckhoff~\cite{kerckhoff}, 
uses analytic properties of the lengths of closed geodesics under earthquakes.
The other, due to Thurston~\cite{thurston-earthquakes}, uses a more geometric
constructions to construct an earthquake from an orientation-preserving
homeomorphism from $S^1$ to itself.

It appears quite likely that those arguments can be extended to provide
other proofs of the ``earthquake theorem'' presented here for hyperbolic
surfaces with geodesic boundary components. The proof given by Thurston, in
particular, might well extend to the case under consideration, however not
in a completely straightforward way since one would have to construct appropriate
homeomorphisms of $S^1$ from two hyperbolic metrics with geodesic boundary and
then understand the boundedness of the earthquake obtained through
Thurston's theorem.

More precisely, given two hyperbolic metrics
$g_1$ and $g_2$ with geodesic boundary on a compact surface with boundary
$\Sigmab$, they define an equivariant self-homeomorphism of the boundary at 
infinity of $\Sigma$, which is a Cantor set in $S^1$. Extending this 
homeomorphism to an equivariant map from $S^1$ to itself can 
be done in many ways. In particular there are $2^k$ such extensions
-- where $k$ is the number of boundary components of $\Sigmab$
corresponding to geodesic boundary components (rather than cusps)
for both $g_1$ and $g_2$ -- obtained by sending each interval in the
complement of the Cantor set to either of its endpoints. It is quite
conceivable that those maps are the boundary values of the earthquakes
considered here. 

In fact this  strategy is not really different from the one we have
considered in this paper. The main point of \cite{thurston-earthquakes} to
construct a left earthquake extending a homeomorphism $\varphi$ of
$S^1_\infty$ is to consider the set $\Ss$ of elements $g$ in $PSL(2,\mr)$ such
that $g\circ\varphi$ is an \emph{extremal} homeomorphism.  The convex hulls in
$\mh^2$ of the fixed points of $g\circ\varphi$ are the strata for the
lamination that provides the earthquake.

Instead, the key point in~\cite{mess}  was to consider the future boundary,
$\partial_+\Kk$, of the convex hull in $AdS_3$ of the
graph of $\varphi$. By means of the product structure of the
boundary of $AdS_3$ two maps were pointed out $M_L,M_R:\partial_+\Kk\rightarrow
P$, where $P\cong\mh^2$ is a fixed spacelike plane. Those maps are
determined by the following requirements:
\begin{enumerate}
\item the restriction on each face is a projective map;
\item ideal points of each face are sent to points on the same left
  (resp. right) leaf.
\end{enumerate}
It turns out that $M_L$ is a left earthquake and $M_R$ is a right earthquake
along the bending lamination of $\partial_+\Kk$ and $M_L\circ M_R^{-1}$
is the earthquake extending $\varphi$.
 
With the AdS language the set $\Ss$ could be identified to the set of points
whose dual plane is a support plane of $\Kk$ touching $\partial_+\Kk$.
Moreover, the intersection of the dual plane with the future boundary is sent
by $M_R$ to the convex core of the fixed points of $g\circ\varphi$.

Let us stress that this relation between these different proofs was already
known by Mess (see the discussion in Section 7 of~\cite{mess}).

\subsubsection*{Other questions.}

Many of the questions which are still open for globally hyperbolic
AdS manifolds (and/or for quasifuchsian hyperbolic manifolds) can also
be considered in the setting of multi-black holes. For instance, Mess~\cite{mess}
asked whether any couple of hyperbolic metrics can be uniquely obtained
as the induced metric on the boundary of the convex core; this might be
true also for hyperbolic metrics with geodesic boundary in the context
of multi-black holes. The corresponding questions for the measured bending
laminations of the boundary of the convex core are also of interest.

\bibliographystyle{amsplain}
\bibliography{../outils/biblio}
\end{document}